\documentclass[11pt]{article}
\usepackage[dvips]{graphicx}
\usepackage{epsf}
\usepackage{amssymb}
\usepackage{pstricks}
\usepackage{pst-node}
\usepackage{lineno}
\usepackage[normalem]{ulem}
\usepackage{mathrsfs}
\usepackage{amsmath}
\usepackage{amsthm}
\usepackage{amsfonts}
\usepackage{eepic}
\usepackage{epic}
\usepackage{color}
\usepackage{amsthm}

\usepackage{lineno}

\newtheorem{thm}{Theorem}[section]

\newtheorem{observation}[thm]{Observation}
\newtheorem{proposition}[thm]{Proposition}
\newtheorem{cor2}[thm]{Corollary}

\newtheorem{remark}{Remark}
\newtheorem{definition}{Definition}
\newcommand{\DTDP}{{\rm DTDP}}
\newcommand{\DT}{{\rm DT}}

\newcommand{\1}{\vspace{0.04cm}}
\newcommand{\2}{\vspace{0.1cm}}

\let\oldenumerate\enumerate
\renewcommand{\enumerate}{
  \oldenumerate
  \setlength{\itemsep}{0.5pt}
  \setlength{\parskip}{0pt}
  \setlength{\parsep}{0pt}
}

\begin{document}

\title{Minimal graphs with disjoint dominating and total dominating sets}
\author{$^{1}$Michael A. Henning\thanks{Research
supported in part by the University of Johannesburg} \, and \, $^{2}$Jerzy Topp\\
\\
$^1$Department of Mathematics and Applied Mathematics\\
University of Johannesburg \\
Auckland Park 2006, South Africa \\
\small {\tt Email: mahenning@uj.ac.za} \\
\\
$^2$State University of Applied Sciences\\ 82-300 Elbl\c{a}g, Poland \\
\small {\tt Email: jtopp@inf.ug.edu.pl}}

\date{}
\maketitle

\begin{abstract}
A graph $G$ is a \DTDP-graph if it has a pair $(D,T)$ of disjoint sets of vertices of $G$ such that $D$ is a dominating set and $T$ is a total dominating set of $G$. Such graphs were studied in a~number of research papers. In this paper we study further properties of \DTDP-graphs and, in particular, we characterize minimal \DTDP-graphs without loops.
\end{abstract}


\section{Introduction}

For notation and graph theory terminology we generally follow~\cite{Henning-Yeo2013}. Let $G=(V_G,E_G)$ be a graph with possible multi-edges and multi-loops. We remark that such a graph is also called a multigraph in the literature. For a vertex $v$ of $G$, its {\em neighborhood\/}, denoted by $N_{G}(v)$, is the set of vertices adjacent to $v$. The {\em closed neighborhood\/} of $v$, denoted by $N_{G}[v]$, is the set $N_{G}(v)\cup \{v\}$. (Observe that if $G$ has a loop incident with $v$, then $N_G(v)=N_G[v]$.) In general, for a subset $X\subseteq V_G$ of vertices, the {\em neighborhood\/} of $X$, denoted by $N_{G}(X)$, is defined to be $\bigcup_{v\in X}N_{G}(v)$, and the {\em closed\/} neighborhood of $X$, denoted by $N_{G}[X]$, is the set $N_{G}(X)\cup X$.

The {\em degree} of a vertex $v$ in $G$, denoted by $d_G(v)$, is the number of edges incident with $v$ plus twice the number of loops incident with $v$.
A vertex of degree one is called a {\em leaf}, and the only neighbor of a~leaf is called its {\em support vertex} (or simply, its {\em support\/}). If a support vertex has at least two leaves as neighbors,  we call it a {\em strong\/} support, otherwise it is a {\em weak} support. The set of leaves, the set of weak supports, the set of strong supports, and the set of all supports of $G$ is denoted by $L_G$, $S'_G$, $S''_G$,  and $S_G$, respectively. If $v$ is a~vertex of $G$, then by $E_G(v)$ and $L_G(v)$ we denote the set of edges and the set of loops incident with $v$ in $G$, respectively.

If $A$ and $B$ are disjoint sets of vertices of $G$, then by  $E_G(A,B)$ is denoted the set of edges in $G$ joining a vertex in $A$ with a~vertex in $B$. For one-element sets we write $E_G(v,B)$,  $E_G(A,u)$, and $E_G(u,v)$ instead of $E_G(\{v\},B)$, $E_G(A,\{u\})$, and $E_G(\{u\},\{v\})$, respectively.

We denote a complete graph, a path, and a cycle on $n$ vertices by $K_n$, $P_n$ and $C_n$, respectively. We note that $C_1$ denoted the cycle of order~$1$ with one loop, and $C_2$ denoted the cycle of order~$2$ with two (repeated) edges. The \emph{corona} $G \circ K_1$ of a graph $G$, also denoted ${\rm cor}(G)$
in the literature, is the graph obtained from $G$ by adding for each vertex $v \in V_G$ a new vertex $v'$ and the edge $vv'$. For an integer $k \ge 1$ we let $[k] = \{1,\ldots,k\}$.

A set of vertices $D \subseteq V_G$ of $G$ is a {\it dominating set\/} if every vertex in $V_G\setminus D$  has a~neighbor in $D$, while $D$ is a {\it total dominating set\/} if every vertex of $G$ has a~neighbor in $D$. We note that a vertex incident with a loop totally dominates all its neighbors (and therefore also itself). A {\em \DT-pair} in a~graph $G$ is a~pair $(D,T)$ of disjoint sets of vertices of $G$ such that $D$ is a dominating set and $T$ is a~total dominating set of $G$. A graph that has a \DT-pair is called a {\em \DTDP-graph} (standing for ``dominating, total dominating, partitionable graph").

Beginning from the classical result of Ore \cite{Ore1962} who was the first to observe that the vertex set of a~graph without isolated vertices can be partitioned into two dominating sets, various graph theoretic properties and parameters of graphs having disjoint dominating sets of different types were studied in a large number of papers.
Properties of \DTDP-graphs (and properties of graphs with disjoint total dominating sets) were extensively studied, for example, in \cite{Broere2004,Delgado2016,Desormeaux2017, Dorfling2005, Henning-Lowenstein-Rautenbach2010-1,Henning-Lowenstein-Rautenbach2010-1+,
Henning-Southey2008,Henning-Southey2009,Henning-Yeo2013,Kiunisala-Jamil,Southey-Henning2011-1}, to mention just a few. In particular, it was proved in [16] that every connected graph with minimum degree at least two and different from $C_5$ is a~\DTDP-graph.  A constructive characterization of all \DTDP-graphs was given in [17].

In this paper we study further properties of \DTDP-graphs and, in particular, we characterize minimal \DTDP-graphs without loops.

\section{Elementary properties of DTDP-graphs}

An immediate consequence of the definition of a \DT-pair is the following observation.

\begin{observation} \label{observ-1} If $(D,T)$ is a \DT-pair in a graph $G$, then every leaf of $G$ belongs to $D$, while every support of $G$ is in $T$, that is, $L_G\subseteq D$ and $S_G\subseteq T$.
\end{observation}

\begin{definition}
{\rm A connected graph $G$ is said to be a \emph{minimal DTDP}-\emph{graph}, if $G$ is a DTDP-graph and no proper spanning subgraph of $G$ is a DTDP-graph. }
\end{definition}

From the definition of a minimal DTDP-graph we immediately have the following observation.

\begin{observation}
\label{observ-2-supergraph}
Every spanning supergraph of a~DTDP-graph is a DTDP-graph, and every DTDP-graph is a spanning supergraph of some~minimal DTDP-graph.
\end{observation}

In view of Observation~\ref{observ-2-supergraph}, minimal DTDP-graphs can be viewed as skeletons of DTDP-graphs, skeletons which can be extended to any DTDP-graph. The following observation determines minimal DTDP-graphs in the class of complete graphs, paths, and cycles.

\begin{observation}
\label{examples-plus}
The following holds. \\ [-24pt]
\begin{enumerate}
\item A complete graph $K_n$ is a DTDP-graph for every $n \ge 3$, but $K_n$ is a minimal DTDP-graph only if $n=3$.
\item A path $P_n$ is a DTDP-graph if and only if $n \in \mathbb{N}\setminus \{1, 2, 3, 5, 6, 9\}$, but $P_n$ is a minimal DTDP-graph only if $n\in \{4, 7, 10, 13\}$.
\item A cycle $C_n$ is a DTDP-graph if $n \ge 3$ and $n\ne 5$, while a~cycle $C_n$ is a~minimal DTDP-graph if and only if $n\in \{3, 6, 9\}$.
\end{enumerate}
\end{observation}

\begin{figure}[h!] \begin{center} \bigskip
{\epsfxsize=4.5in \epsffile{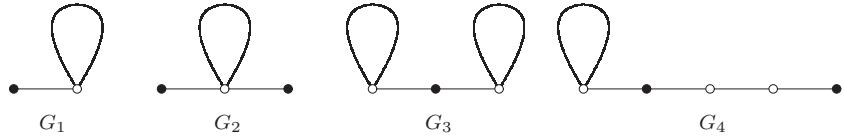}}

\vspace{-2mm}
\caption{Small minimal $D\!T\!D\!P$-graphs with loops} \label{minimal-graphs-with-loops}
\end{center}\end{figure}

Examples of small minimal DTDP-graphs with loops are given in Fig.~\ref{minimal-graphs-with-loops}.

Let $G$ be a graph and let $v$ be a support vertex in $G$. If $G'$ is obtained from $G$ by adding a new leaf, say $v'$, adjacent to $v$, then $(D,T)$ is a \DT-pair in $G$ if and only if $(D\cup\{v'\},T)$ is a \DT-pair in $G'$. This yields the following observation, which shows that the addition of new leaves adjacent to an existing support vertex of a graph preserves the property of being a (minimal)  DTDP-graph.

\begin{observation}
\label{observ-S2(H)-stron-weak-supports}
Let $v$ be a support vertex in a~graph $G$. If $G'$ is a graph obtained from $G$ by adding a new leaf, say $v'$, adjacent to $v$, then $G$ is a~{\rm (}minimal\,{\rm )} DTDP-graph if and only if $G'$ is a~{\rm (}minimal\,{\rm )}  DTDP-graph.
\end{observation}

The following result determines the minimal DTDP-graphs in the class of coronas of graphs.

\begin{proposition} 
\label{observ-1-dla korony}
If $H$ is a connected graph of size at least one, then the corona $H \circ K_1$ is a~DTDP-graph. Moreover, the corona $H\circ K_1$ is a minimal DTDP-graph if and only if $H$ is a star or $C_1$.
\end{proposition}
\begin{proof} Assume that $H$ is a connected graph of size at least one. The corona $H \circ K_1$ of $H$ is a~DTDP-graph, as the sets $D = L_{H \circ K_1}$ and $T = V_{H \circ K_1} \setminus D$ form a \DT-pair in $H \circ K_1$. Assume now that $H$ is a star or $C_1$. Let $e$ be an arbitrary edge or a~loop (in the case when $H = C_1$) in $H \circ K_1$. We note that $K_1$ or $K_2$ is a component of $(H \circ K_1)-e$, and therefore $(H\circ K_1)-e$ is not a DTDP-graph, implying that $H \circ K_1$ is a minimal DTDP-graph.

Finally assume that $H$ is a connected graph of size at least one and $H \circ K_1$ is a~minimal DTDP-graph. We claim that $H$ is a star or $C_1$. If $H$ has order~$1$, then $H \circ K_1$ is a minimal DTDP-graph if and only if $H=C_1$. Thus, we may assume that $H$ is a graph of order at least~$2$, for otherwise the desired result follows. Let $e$ be an arbitrary edge or a loop in $H$. If $e$ is not a pendant edge in $H$, then $H-e$ is without isolated vertices, implying that $(H \circ K_1) - e = (H-e) \circ K_1$ is a DTDP-graph, contradicting the minimality of $H \circ K_1$. Thus, $e$ is a pendant edge in $H$. Since $e$ is an arbitrary edge in $H$, the graph $H$ is a star. \end{proof}

Let $k$ be a nonnegative integer, and let ${\cal S}_k$ denote the family of all rooted trees in which all leaves are at distance $k$ from the root.

\begin{proposition} \label{observ-sniezynki}
No tree belonging to ${\cal S}_0\cup {\cal S}_1\cup {\cal S}_2$ is a DTDP-graph, while every tree belonging to the family ${\cal S}_3$ is a minimal DTDP-graph.
\end{proposition}
\begin{proof} It is evident (and follows from Observation \ref{observ-1}) that no tree belonging to ${\cal S}_0\cup {\cal S}_1\cup {\cal S}_2$ is a~DTDP-graph. Now let $T$ be a tree belonging to ${\cal S}_3$. The root of $T$, together with the vertices on level~$3$ (at distance~$3$ from $v$ in $T$) form a dominating set of $T$, while the vertices on level~$1$ and~$2$ (at distance~$1$ and~$2$, respectively, from~$r$ in $T$) form a total dominating set of $T$. This proves that $T$ is a  DTDP-graph. Its minimality follows from the fact that every proper spanning subgraph of $T$ contains a component belonging to ${\cal S}_0\cup {\cal S}_1\cup {\cal S}_2$, and is therefore not a DTDP-graph.
\end{proof}

\section{2-Subdivision graphs of a graph}%

Let $H=(V_H,E_H)$ be a graph  with possible multi-edges and multi-loops. By $\varphi_H$ we denote a~function from $E_H$ to $2^{V_H}$ that associates with each $e \in E_H$, the set $\varphi_H(e)$ of vertices incident with $e$. Let $X_2$ be a set of $2$-element subsets of an arbitrary set (disjoint with $V_H\cup E_H$), and let $\xi \colon E_H \to X_2$ be a function such that $\xi(e) \cap \xi(f) = \emptyset$ if $e$ and $f$ are distinct elements of $E_H$. If $e \in E_H$ and $\varphi_H(e) = \{u,v\}$ ($\varphi_H(e) = \{v\}$, respectively), then we write $\xi(e) = \{u_e,v_e\}$ ($\xi(e) = \{v^1_e,v^2_e\}$, respectively). Let $S_2(H)$ denote the graph obtained from $H$ by inserting two new vertices into each edge and each loop of $H$. Formally, the graph $S_2(H)$ has vertex set
\[
V_{S_2(H)}= V_H\cup  \bigcup_{e\in E_H}\xi(e)
\]
and edge set $E_{S_2(H)} = E_1 \cup E_2$, where
\[
\begin{array}{lcl}
E_1 & = & \displaystyle{ \bigcup_{e\in E_H}\{xy\colon \xi(e)=\{x,y\}\}, \, \mbox{ and} } \2 \\
E_2 & = & \displaystyle{ \bigcup_{v\in V_H} \big(\{vv_e \colon e\in E_H(v)\} \cup \{vv^1_e, vv^2_e \colon e \in L_H(v)\}\big). } \\
\end{array}
\]

In such a graph $S_2(H)$, we let ${\cal P} = \{{\cal P}(v) \colon v\in V_H\}$ be a family in which ${\cal P}(v)$ is a~partition of the set $N_{S_2(H)}(v)$ for each $v\in V_H \subset V_{S_2(H)}$. Further, we let $S_2(H,{\cal P})$ denote the graph (possibly with multi-loops and multi-edges) with vertex set
\[
V_{S_2(H,{\cal P})} = V_H \cup \bigcup_{v \in V_H}(\{v\}\times {\cal P}(v))
\]
and edge set $E_{S_2(H,{\cal P})}$ defined by the neighborhoods of the vertices of $S_2(H,{\cal P})$, where
\[
N_{S_2(H,{\cal P})}(v) = \{(v,A)\colon A\in {\cal P}(v)\}
\]
if $v \in V_H \subseteq V_{S_2(H,{\cal P})}$, and
\[
N_{S_2(H,{\cal P})}((v,A)) = \{v\}\cup  \bigcup_{u\in N_H(v)} \{(u,B) \colon B\in {\cal P}(u)\, \mbox{and}\, N_{S_2(H)}(A)\cap B\ne \emptyset\}
\]
if $(v,A) \in \bigcup_{v\in V_H}(\{v\}\times {\cal P}(v))$. In the last case it is necessary to add that if $u$ and $v$ are adjacent vertices in $H$, $A\in {\cal P}(v)$, and $B\in {\cal P}(u)$, then the number of edges joining $(v,A)$ and $(v,B)$ in $S_2(H,{\cal P})$ is equal to the cardinality of the set $\{e\in E_H \colon \varphi_H(e)=\{u,v\},\,\, v_e \in A,\,\, \mbox{and}\,\, u_e\in B\}$. Similarly, we can determine the number of edges between vertices $(v,A)$ and $(v,B)$ (and the number of loops incident with $(v,A)$) if $A, B\in {\cal P}(v)$. Intuitively, $S_2(H,{\cal P})$ is the graph obtained from $S_2(H)$ by contracting the vertices of the set $A\in {\cal P}(v)$ into a new vertex $(v,A)$, which becomes adjacent to all former neighbors of the vertices belonging to $A$ (for every $A\in {\cal P}(v)$ and every $v\in V_H\subseteq V_{S_2(H)}$). It is evident that $S_2(H,{\cal P})$ is isomorphic to $S_2(H)$ if ${\cal P}= \{{\cal P}(v)\colon v\in V_H\}$ and ${\cal P}(v)=\{\{x\}\colon x\in N_H(v)\}$ for each $v\in V_H$.

\begin{figure}[h!] \begin{center} \bigskip
{\epsfxsize=5.99in \epsffile{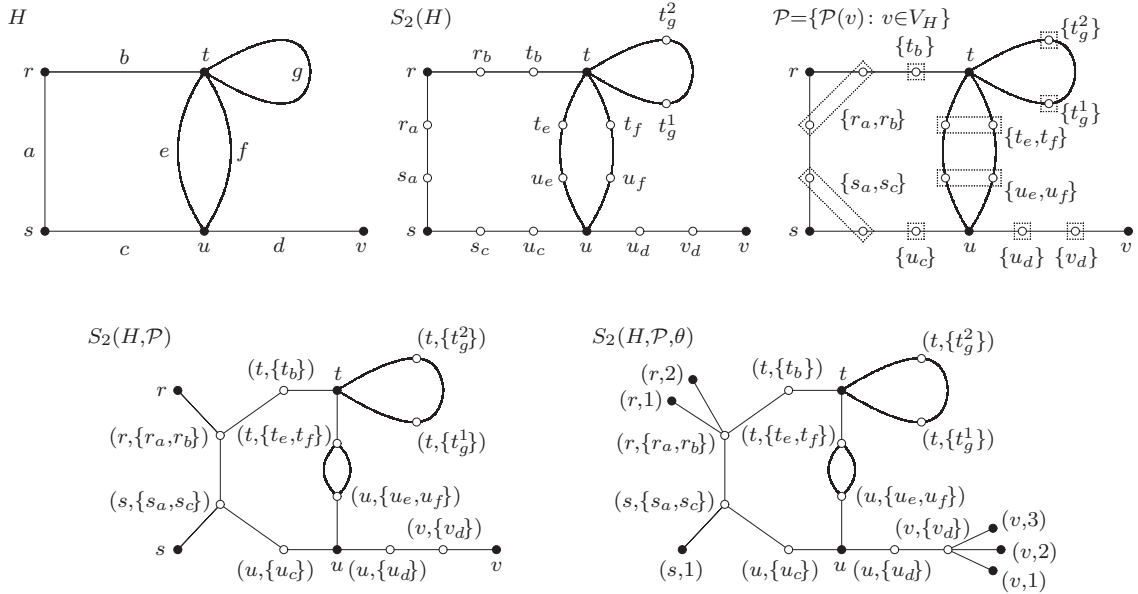}}

\vspace{-2mm}
\caption{2-subdivision graphs $S_2(H)$, $S_2(H,{\cal P})$, and $S_2(H,{\cal P},\theta)$ of a graph $H$} \label{2-subdivision-graph-new}
\end{center}\end{figure}

For a positive function $\theta \colon L_{S_2(H,{\cal P})} \to \mathbb{N}$, we let
\[
L_\theta \colon L_{S_2(H,{\cal P})} \to L_{S_2(H,{\cal P})} \times \mathbb{N}
\]
be a~function such that $L_\theta(x)=\{(x,i) \colon i \in [\theta(x)] \}$ for $x \in L_{S_2(H,{\cal P})}$. Finally, let ${S_2(H,{\cal P},\theta)}$ denote the graph obtained from ${S_2(H,{\cal P})}$ by replacing each leaf $x$ of ${S_2(H,{\cal P})}$ by its copies $(x,1),\ldots, (x,\theta(x))$ (adjacent to the former neighbor of $x$ in ${S_2(H,{\cal P})}$), see, for example, Fig. \ref{2-subdivision-graph-new}.  Certainly, ${S_2(H,{\cal P},\theta)}$ is isomorphic to ${S_2(H,{\cal P})}$ if $\theta(x)=1$ for every leaf $x$ of ${S_2(H,{\cal P})}$. The three graphs $S_2(H)$, ${S_2(H,{\cal P})}$, and ${S_2(H,{\cal P},\theta)}$ are said to be {\em $2$-subdivision graphs} of $H$ (for a family ${\cal P}=\{{\cal P}(v)\colon v\in V_H\}$ of partitions ${\cal P}(v)$ of neighborhoods $N_{S_2(H)}(v)$ where $v\in V_H\subset V_{S_2(H)}$), and for a~positive function $\theta \colon L_{S_2(H,{\cal P})}\to \mathbb{N}$$)$. Let $V_{S_2(H)}^o$, $V_{S_2(H,{\cal P})}^o$, and $V_{S_2(H,{\cal P},\theta)}^o$ be sets of vertices such that \[
\begin{array}{lcl}
V_{S_2(H)}^o & = & V_{S_2(H)}\setminus  V_{S_2(H)}^n, \1 \\
V_{S_2(H,{\cal P})}^o & = & V_{S_2(H,{\cal P})}\setminus V_{S_2(H,{\cal P})}^n, \1 \\
V_{S_2(H,{\cal P},\theta)}^o & = & V_{S_2(H,{\cal P},\theta)}\setminus V_{S_2(H,{\cal P},\theta)}^o, \2
\end{array}
\]
where $V_{S_2(H)}^n=\bigcup_{e\in E_H}\xi(e)$, and $V_{S_2(H,{\cal P})}^n= V_{S_2(H,{\cal P},\theta)}^n = \bigcup_{v\in V_H}(\{v\}\times {\cal P}(v))$.

As a consequence of the definition of the $2$-subdivision graphs $S_2(H)$ and $S_2(H,{\cal P})$, we have the following observation.

\begin{observation} 
\label{pierwsze-wlasnosci-S2(H)}
If $H$ is a graph with no isolated vertex and ${\cal P}=\{{\cal P}(v)\colon v\in V_H\}$ is a family in which ${\cal P}(v)$ is a partition of the neighborhood $N_{S_2(H)}(v)$ for each $v \in  V_H \subset V_{S_2(H)}$, then the following statements hold. \\ [-22pt]
\begin{enumerate}
\item[$(1)$] If $v\in V_H$, then $N_{S_2(H)}(v)$ and $N_{S_2(H,{\cal P})}(v)$ are nonempty subsets of $V_{S_2(H)}\setminus V_H$ and  $V_{S_2(H,{\cal P})}\setminus V_H$, respectively.
\item[$(2)$] If $v\in V_H$, then $d_{S_2(H)}(v) = d_H(v)$ and $d_{S_2(H,{\cal P})}(v)= |{\cal P}(v)|$.
\item[$(3)$] If $x\in V_{S_2(H)}\setminus V_H$, then $|N_{S_2(H)}(x)\cap V_H|= |N_{S_2(H)}(x)\cap (V_{S_2(H)}\setminus V_H)|=1$.
\item[$(4)$] If $v\in V_H$ and $A\in {\cal P}(v)$, then $|N_{S_2(H,{\cal P})}((v,A))\cap (V_{S_2(H,{\cal P})}\setminus  V_H)|=|A|$ and $|N_{S_2(H,{\cal P})}((v,A))\cap V_H|=1$.
\end{enumerate}
\end{observation}

We prove next the following result.

\begin{proposition} 
\label{observ-S2(H)-jest-DPDP-grafem}
If no vertex of a graph $H$ is an isolated vertex, then the $2$-subdivision graph $S_2(H,{\cal P},\theta)$ is a DTDP-graph for every family ${\cal P}=\{ {\cal P}(v)\colon v\in V_H\}$ of partitions ${\cal P}(v)$ of neighborhoods $N_{S_2(H)}(v)$  where $v\in V_H$, and for every positive function $\theta \colon L_{S_2(H,{\cal P})}\to \mathbb{N}$. In addition, $(V_{S_2(H,{\cal P},\theta)}^o,V_{S_2(H,{\cal P},\theta)}^n)$ is a  \DT-pair in $S_2(H,{\cal P},\theta)$.
\end{proposition}
\begin{proof} It follows from the statements (1) and (4) of Observation \ref{pierwsze-wlasnosci-S2(H)} that $V_{S_2(H,{\cal P})}\setminus V_H$ is a total dominating set of $S_2(H,{\cal P})$. The second part of Observation \ref{pierwsze-wlasnosci-S2(H)}\,(4) also proves that $V_H$ is a~dominating set of $S_2(H,{\cal P})$. Therefore, $S_2(H,{\cal P})$ is a DTDP-graph. Consequently, $S_2(H,{\cal P},\theta)$ is a~DTDP-graph if $\theta(x)=1$ for every $x\in L_{S_2(H,{\cal P})}$ (since $S_2(H,{\cal P})$ and $S_2(H,{\cal P},\theta)$ are isomorphic if $\theta(x)=1$ for every $x\in L_{S_2(H,{\cal P})}$). From this and from Observation \ref{observ-S2(H)-stron-weak-supports}, it follows that $S_2(H,{\cal P},\theta)$ is a DTDP-graph for every positive function $\theta \colon L_{S_2(H,{\cal P})}\to \mathbb{N}$. In addition, it is obvious that  $(V_{S_2(H)}^o,V_{S_2(H)}^n)$,
$(V_{S_2(H,{\cal P})}^o,V_{S_2(H,{\cal P})}^n)$, and $(V_{S_2(H,{\cal P},\theta)}^o, V_{S_2(H,{\cal P},\theta)}^n)$ are \DT-pairs in $S_2(H)$, $S_2(H,{\cal P})$, and $S_2(H,{\cal P},\theta)$, respectively. \end{proof}

\begin{figure}[h!] \begin{center} \bigskip
{\epsfxsize=5.0in \epsffile{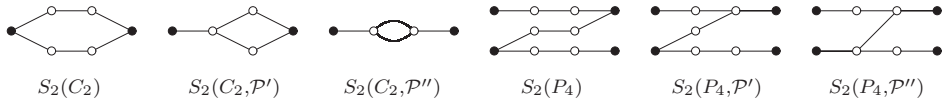}}

\vspace{-2mm}
\caption{Graphs $S_2(C_2)$, $S_2(C_2,{\cal P}')$, $S_2(C_2,{\cal P}'')$,  $S_2(P_4)$, $S_2(P_4,{\cal P}')$, and $S_2(P_4,{\cal P}'')$} \label{examples}
\end{center}\end{figure}

\begin{figure}[h!] \begin{center} \bigskip
{\epsfxsize=5in \epsffile{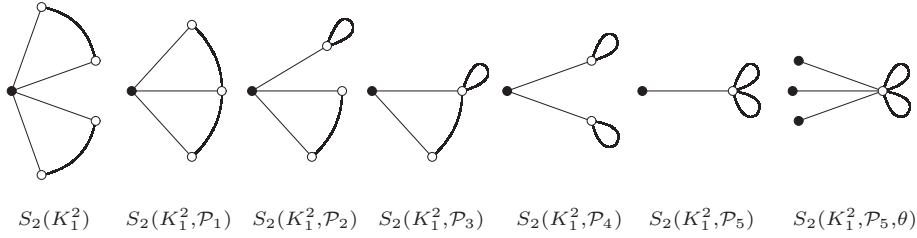}}

\vspace{-2mm}
\caption{Possible 2-subdivision graphs $S_2(K_1^2,{\cal P},\theta)$ of $K_1^2$} \label{examples-K-1-m}
\end{center}\end{figure}

By Proposition~\ref{observ-S2(H)-jest-DPDP-grafem} every $2$-subdivision graph $S_2(H,{\cal P},\theta)$ is a DTDP-graph. Simple examples presented in Fig. \ref{examples} and
\ref{examples-K-1-m} illustrate the fact that if $S_2(H,{\cal P},\theta)$ is a minimal DTDP-graph depends on the family of partitions ${\cal P}$. In Fig. \ref{examples-K-1-m} we present examples of possible $2$-subdivision graphs $S_2(K_1^2,{\cal P},\theta)$ of $K_1^2$, where $K_1^s$ denotes a graph of order 1 and size $s$. It is easy to check that of all these  graphs in Fig. \ref{examples} and \ref{examples-K-1-m} only $S_2(C_2)$, $S_2(P_4)$, and $S_2(K_1^2,{\cal P}_4)$ are minimal  DTDP-graphs. In the next two observations we study these relations more precisely.

\begin{proposition} 
\label{o sklejaniu krawedzi} Let $H$ be a graph without isolated vertices and let ${\cal P}=\{{\cal P}(v)\colon v\in V_H\}$ be a~family in which ${\cal P}(v)$ is a~partition of the set $N_{S_2(H)}(v)$ for $v\in V_H$. Then the $2$-subdivision graph $S_2(H,{\cal P})$ $($and $S_2(H,{\cal P},\theta)$ for every positive function $\theta \colon L_{S_2(H,{\cal P})}\to \mathbb{N}$$)$ is a DTDP-graph but not a~minimal DTDP-graph if there exists a vertex $v\in V_H$ and a~set $A\in {\cal P}(v)$ such that $A$ contains at least two elements belonging to different edges or loops and at least one of them is not a pendant edge of $H$.
\end{proposition}
\begin{proof} By Proposition~\ref{observ-S2(H)-jest-DPDP-grafem}, $S_2(H,{\cal P})$
is a DTDP-graph and the pair $(D,T)= (V_{S_2(H,{\cal P})}^o,$ $V_{S_2(H,{\cal P})}^n)$ is a~\DT-pair in $S_2(H,{\cal P})$. It suffices to observe that some proper spanning subgraph of $S_2(H,{\cal P})$ is a DTDP-graph. We consider two possible cases.

{\it Case 1}.
Assume that $e$ is a non-pendant edge in $H$, say $\varphi_H(e)=\{v,u\}$, and let
$f$ be an edge (or a loop) incident with $v$ and such that $v_e$ and $v_f$ ($v_e$ and $v_f^1$ or $v_e$ and $v_f^2$, if $f$ is a loop) belong to the same set $A$ which is an element of the partition ${\cal P}(v)$. Let $B$ be the only set belonging to ${\cal P}(u)$ that contains $u_e$. Now from the properties of the \DT-pair $(D,T)$ in $S_2(H,{\cal P})$ it follows that if $\{u_e\} \varsubsetneq B$, then $(D,T)$ is a~\DT-pair in the proper spanning subgraph $S_2(H,{\cal P})\setminus (v,A)(u,B)$ of $S_2(H,{\cal P})$. Similarly, if $B=\{u_e\}$, then $(D\cup\{(u,\{u_e\})\}, T\setminus \{(u,\{u_e\})\})$ is a \DT-pair in the proper spanning subgraph $S_2(H,{\cal P})- u(u,\{u_e\})$ of $S_2(H,{\cal P})$.

{\it Case 2}. Assume now that $e$ is a loop incident with a vertex $v$ in $H$, and $f$ is a~pendant edge or a loop incident with $v$ in $H$ ($f\ne e$) and such that $\{v_e^1, v_e^2\}\cap A\ne  \emptyset$ and $v_f\in A$ (or $\{v_e^1, v_e^2\}\cap A\ne  \emptyset$ and $\{v_f^1, v_f^2\}\cap A\ne  \emptyset$, if $f$ is a loop) for some $A\in {\cal P}(v)$. If $f$ is a~pendant edge, and, without loss of generality, if $v_e^1, v_f \in A$, then we consider three subcases: $(i)$ $v_e^2 \in A$; $(ii)$ $v_e^2 \not\in A$ and $\{v_e^2\}\in {\cal P}(v)$; $(iii)$ $v_e^2 \not\in A$ and $\{v_e^2\}\varsubsetneq B\in {\cal P}(v)$. In the first case $(D,T)$ is a~\DT-pair in the spanning subgraph obtained from $S_2(H,{\cal P})$ by removing one loop incident with the vertex $(v,A)$. In the second case the pair $(D\cup\{(v,\{v_e^2\})\}, T\setminus \{(v,\{v_e^2\})\})$ is a \DT-pair in $S_2(H,{\cal P})- v(v,\{v_e^2\})$. In the third case $(D,T)$ is a~\DT-pair in $S_2(H,{\cal P})-(v,A)(u,B)$.

Let us consider graphs $H_v=H[\{v\}]$, $S_2(H_v)$, and $S_2(H_v, {\cal P}_v)$, where ${\cal P}_v$ is the partition of $N_{S_2(H_v)}(v)$ defined by ${\cal P}(v)$, i.e.,
${\cal P}_v= \{A\cap N_{S_2(H_v)}(v)\colon A\in {\cal P}(v),\,\, A\cap N_{S_2(H_v)}(v)\ne  \emptyset\}$. It is obvious that $S_2(H_v, {\cal P}_v)$ is an induced subgraph of $S_2(H, {\cal P})$, and $(\{v\}, V_{S_2(H_v, {\cal P}_v)}\setminus \{v\})$ is a \DT-pair in $S_2(H_v, {\cal P}_v)$. It remains to prove that some proper spanning subgraph of $S_2(H_v, {\cal P}_v)$ has a \DT-pair $(D_v,T_v)$ such that $v \in D_v$. This is obvious if $S_2(H_v, {\cal P}_v)$ contains parallel edges or a~loop incident with a vertex of degree at least 2. Thus assume that $S_2(H_v, {\cal P}_v)$ contains neither a loop incident with a vertex of degree at least 2 nor parallel edges. We may also assume that no two mutually adjacent vertices of degree 2 are adjacent to $v$. Let $u\in N_{S_2(H_v)}(v)$ be a vertex of minimum degree in $S_2(H_v, {\cal P}_v)$. Certainly, $d_{S_2(H_v, {\cal P}_v)}(u)\ge 2$. Let $w$ be a~vertex belonging to
$N_{S_2(H_v, {\cal P}_v)}(u)\setminus \{v\}$. The choice of $u$ and the assumption that no two mutually adjacent vertices of degree 2 are adjacent to $v$ imply that $w$ is of degree at least 3 in $S_2(H_v, {\cal P}_v)$. Consequently, $w$ has a neighbor in $V_{S_2(H_v, {\cal P}_v)}\setminus \{v,u\}$. This implies that $(\{v,u\}, V_{S_2(H_v, {\cal P}_v)}\setminus \{v,u\})$ is a \DT-pair in $S_2(H_v, {\cal P}_v)-uv$, and completes the proof.
\end{proof}

\begin{proposition} 
\label{o sklejaniu krawedzi koncowych} Let $H$ be a~graph without isolated vertices, let ${\cal P}=\{{\cal P}(v)\colon v\in V_H\}$ be a family in which ${\cal P}(v)$ is a partition of the set $N_{S_2(H)}(v)$ for every $v\in V_H$. If ${\cal P}$ is such that for every $v\in V_H\setminus L_H$ and every non-pendant edge $e$ incident with $v$, the singleton $\{v_e\}$ is an element of ${\cal P}(v)$, then both $2$-subdivision graphs $S_2(H)$ and $S_2(H,{\cal P})$ are minimal DTDP-graphs  or neither of them is a minimal DTDP-graph.
\end{proposition}
\begin{proof} For ease of observation, we assume that $H$ has only one support vertex, say $v$. Let $u^1, \ldots, u^k$ be the leaves adjacent to $v$ in $H$. Let $H'$ be the subgraph of $H$ induced by the vertices $v, u^1, \ldots, u^k$. It follows from the properties of ${\cal P}$ that the $2$-subdivision graph $S_2(H,{\cal P})$ results from  $S_2(H)$ replacing the tree $S_2(H')$ rooted at $v$ by the tree $S_2(H',{\cal P}')$ rooted at $v$ and defined for the family ${\cal P}'=\{{\cal P}'(x)\colon x\in V_{H'}\}$
in which ${\cal P}'(v)=\{A\in {\cal P}\colon A\subseteq \{v_{vu^1},\ldots, v_{vu^k}\}\}$ and ${\cal P}'(u^i)= {\cal P}(u^i)= \{u^i_{vu^i}\}$ (for $i\in [k]$).
Now, the fact that both $S_2(H')$ and $S_2(H',{\cal P}')$ are minimal DTDP-graphs (by Proposition~\ref{observ-sniezynki}) implies that both $S_2(H)$ and $S_2(H,{\cal P})$ are minimal DTDP-graphs  or neither of them is a minimal DTDP-graph.
\end{proof}

If $e$ is a pendant edge in $H$ and $\varphi_H(e)=\{v,u\}$, where $v$ is a support vertex of degree at least~$2$ and $u$ is a leaf, then the edge $vv_e$ in $S_2(H)$ is called a {\em far part of the pendant edge $e$} in $H$. If $e$ is a loop incident with a vertex $v$ in $H$, then the edges $vv_e^1$ and $vv_e^2$ in $S_2(H)$ are said to be {\em twin parts of the loop $e$} in $H$. It follows from Propositions \ref{o sklejaniu krawedzi} and \ref{o sklejaniu krawedzi koncowych}  that if $S_2(H,{\cal P},\theta)$ is a minimal DTDP-graph, then ${\cal P}$ can only contract far parts of adjacent pendant edges in $H$ or twin parts of a loop in $H$.

\begin{cor2} If $H$ is a connected graph of size at least~$2$, then the $2$-subdivision graphs  $S_2(H,{\cal P})$ and $S_2(H,{\cal P},\theta)$ are minimal DTDP-graphs
$($for every family ${\cal P}=\{{\cal P}(v)\colon v\in V_H\}$ in which ${\cal P}(v)$ is a~partition of the set $N_{S_2(H)}(v)$ for each $v\in V_H$ and for every positive function $\theta \colon L_{S_2(H,{\cal P})}\to \mathbb{N}$$)$ if and only if $H$ is a star. \end{cor2}

Our aim is to recognize graphs which are present in non-minimal DTDP-graphs. For this purpose, let ${\cal F}$ be the family of all graphs defined as follows:
\begin{enumerate}
\item[$(1)$] We start with a rooted tree, say  $T$, in which $d_T(x)\le 2$ for every vertex $x\in V_T\setminus \{r\}$ (where $r$ is a root of $T$) and form $2$-subdivision graphs $S_2(T)$ and $S_2(T,{\cal P})$ for the family ${\cal P}=\{{\cal P}(x)\colon x\in V_T\}$ in which ${\cal P}(x)$ is a partition of the set $N_{S_2(T)}(x)$, and where ${\cal P}(r)=\{N_{S_2(T)}(r)\}$, while ${\cal P}(x)= \{\{y\}\colon y\in N_{S_2(T)}(x)\}$ if $x\in V_T\setminus \{r\}$.
\item[$(2)$] From the choice of $\cal P$ it follows that $r$ is a leaf in $S_2(T,{\cal P})$ (and $(r,N_{S_2(T)}(r))$ is the only neighbor of $r$ in $S_2(T,{\cal P})$. Finally, if $\theta \colon L_{S_2(T,{\cal P})}\to \mathbb{N}$ is a function such that $\theta(r)$ is a positive integer and $\theta(x)=1$ for every $x\in L_{S_2(T,{\cal P})}\setminus \{r\}$, then we form the $2$-subdivision graph $S_2(T,{\cal P},\theta)$ from $S_2(T,{\cal P})$ replacing the leaf $r$ by its copies $(r,1),\ldots, (r,\theta(r))$ (adjacent to $(r,N_{S_2(T)}(r))$ in ${S_2(H,{\cal P},\theta)}$).
\end{enumerate}

\begin{figure}[h!] \begin{center} \bigskip
{\epsfxsize=5.25in \epsffile{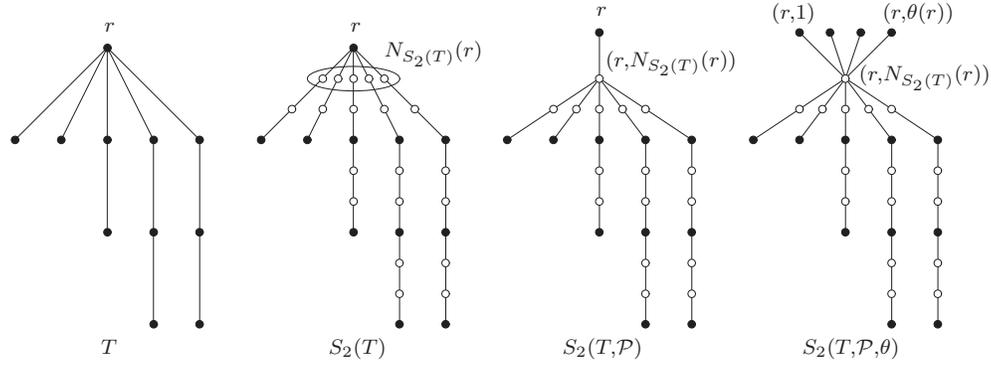}}

\vspace{-2mm}
\caption{A graph $S_2(T,{\cal P},\theta)$ belongs to the family ${\cal F}$} \label{Remark-przed-a-good-subgraph}
\end{center}\end{figure}

To illustrate this definition, consider graphs drawn in Fig. \ref{Remark-przed-a-good-subgraph}. We remark that the graph $S_2(T,{\cal P},\theta)$ belonging to the family ${\cal F}$ can be obtained from the graph $S_2(T)$ by removing the leaves different from the root $r$ of $T$, and then adding a positive integer of new leaves adjacent to the vertex $r$. From this remark and from Proposition~\ref{observ-S2(H)-jest-DPDP-grafem} we have the following corollary, which we shall use  in the proof of Theorem \ref{zakazane-w-minimalnych-DTDP}.

\begin{cor2} \label{Corollary family F} Every tree belonging to the family $\cal F$ is a DTDP-tree, that is, if $T$ is a tree rooted at vertex $r$, and $d_T(x)\le 2$ for every vertex $x\in V_T\setminus \{r\}$,  $1\le |N_T(r)\cap L_T|\le d_T(r)-1$, and $d_T(x,r)\equiv 2\,(\!\!\!\!\mod 3)$ for every $x\in L_T\setminus N_T(r)$, then $T$ is a  DTDP-tree. In particular, if $T$ is a~wounded spider rooted at vertex $r$, that is, if $1\le |N_T(r)\cap L_T|\le d_T(r)-1$ and $d_T(x,r)= 2$ for every $x\in L_T\setminus N_T(r)$, then $T$ is a~minimal DTDP-tree.\end{cor2}

\vspace{-4mm}

\section{Good subgraphs of a graph}

In this section, we define a good subgraph of a graph. Let $Q$ be a subgraph without isolated vertices of a graph $H$, and let $E_Q^-$ denote the set of edges belonging to $E_H \setminus E_Q$ that are incident with a~vertex of $Q$. Let $E$ be a set such that $E_Q^-\subseteq E \subseteq E_H \setminus E_Q$, and let $A_E$ denote a set of arcs obtained by assigning an orientation for each edge in $E$. Then by  $H(A_E)$ we denote the partially oriented graph obtained from $H$ by replacing the edges in $E$ by the arcs belonging to $A_E$. If $e\in E$, then by $e_A$ we denote the only arc in $A_H$ that corresponds to $e$. By $H_0$ we denote the subgraph of $H(A_E)$ induced by the vertices that are not the initial vertex of an arc belonging to $A_E$, i.e., by the set $\{v\in V_H \colon d^+_{H(A_E)}(v) = 0\}$.

We say that $Q$ is a~\emph{good subgraph of $H$} if there exist a set of edges $E$ (where $E_Q^-\subseteq E \subseteq E_H \setminus E_Q$) and a set of arcs $A_E$ such that in the resulting graph $H(A_E)$, which we simply denote by $H$ for notational convenience, the arcs in $A_E$ form a~family ${\cal F}= \{F_v \colon v\in V_Q\}$ of arc disjoint digraphs $F_v$ indexed by the vertices of $Q$ and such that the following holds.
\\ [-20pt]
\begin{enumerate}
\item[{\rm (1)}] For every $v$ in $Q$, the digraph $F_v$ is the union of a family, say ${\cal P}_v$, of arc disjoint oriented paths that begin at $v$.
\item[{\rm (2)}]  If $u\in V_H\setminus V_Q$, then $d^+_H(u)\le 1$.
\item[{\rm (3)}] If $u\in V_H$, then $d^-_H(u)<d_H(u)$.
\item[{\rm (4)}] If $x \in V_{F_v}\cap V_{F_u}$ and $v\ne u$, then $d_{F_v}^+(x)=0$ or $d_{F_u}^+(x)=0$.
\end{enumerate}

\begin{figure}[h!] \begin{center} \bigskip
{\epsfxsize=5.5in \epsffile{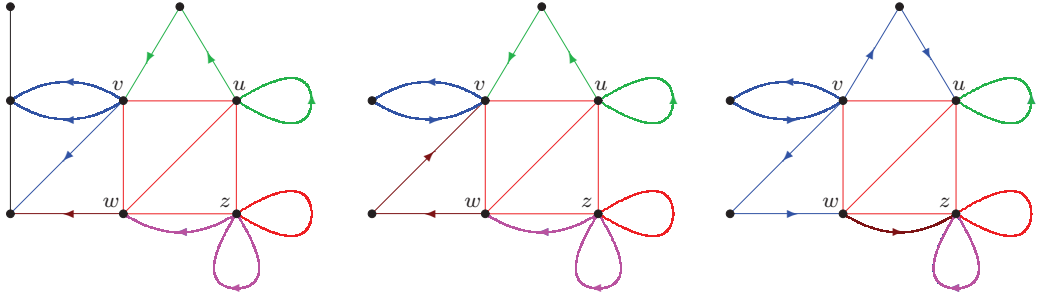}}

\vspace{-2mm}
\caption{A red subgraph is a good subgraph in the host graph} \label{yyy-graphs}
\end{center}\end{figure}

One example of a good subgraph $Q$ (in different graphs) is presented in Fig. \ref{yyy-graphs}. For clarity, the edges of $Q$ are drawn in red, the digraphs $F_v$, $F_u$, $F_w$, and $F_z$ are drawn in blue, green, brown, and purple, respectively,  orientations of arcs are represented by arrows, and all other edges (here only two) are black.

From the  definition of a~good subgraph we immediately have the following observation.

\begin{observation} \label{leaves-supports} Neither a~leaf nor a~support vertex of a graph $H$ belongs to a~good subgraph in $H$. \end{observation}

Observation \ref{leaves-supports} also implies that not every graph has a good subgraph. In particular, a corona graph (that is, a graph in which each vertex is a leaf or it is adjacent to exactly one leaf) has no good subgraph.  On the other hand, if $Q$ is a graph with no isolated vertex, and $H$ is the graph obtained from $Q$ by attaching at least one pendant edge to each vertex of $Q$ and thereafter subdividing these new edges, then $Q$ is a~good subgraph in $H$. This proves that every graph without isolated vertices can be a good subgraph of some graph.

\begin{figure}[h!] \begin{center} \bigskip
{\epsfxsize=4.0in \epsffile{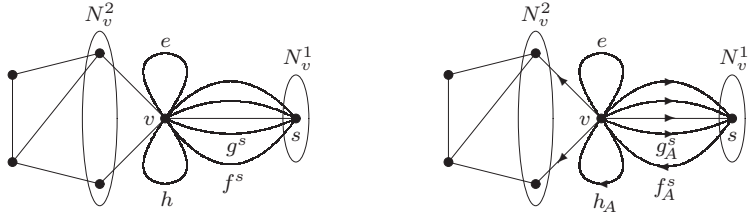}}

\vspace{-2mm}
\caption{An example to Observation \ref{remark1}} \label{Rys-remark 1}
\end{center}\end{figure}

\begin{proposition}
\label{remark1} If $e$ is a~loop incident with a vertex $v$ in a connected graph $H$, then the subgraph $H_e$ of $H$ generated by $e$ is a good subgraph in $H$ if and only if $H\ne  C_1$ and $v$ is not adjacent to a~pendant edge in $H$.  \end{proposition}
\begin{proof}
If $H_e$ is a good subgraph in $H$, then, by Observation \ref{leaves-supports}, the only vertex $v$ of $H_e$ cannot be a support vertex in $H$.

Assume now that the vertex $v$ is not a support vertex in $H$. Then $v$ is neither a~support vertex nor a leaf in $H$ (as $e$ is a loop incident with $v$). If $N_H(v) =\{v\}$, then $H=K_1^s$ ($s\ge 2$) and, certainly, $H_e=K_1^1$ is a good subgraph in $H$. Thus assume that $N_H(v)\ne \{v\}$. In this case, the set $N_H(v)\setminus \{v\}$ is nonempty and it consists of two disjoint subsets $N_v^1$ and $N_v^2$, where $N_v^1= \{x\in N_H(v)\setminus \{v\}\colon N_H(x)= \{v\}\}$ and $N_v^2= \{x\in N_H(v)\setminus \{v\} \colon \{v\} \varsubsetneq N_H(x)\}$. Consequently, the set $E_{H_e}^-$ of edges or loops belonging to $E_H\setminus E_{H_e}= E_H\setminus \{e\}$ that are incident with $v$, consists of three disjoint subsets $E_v^l$,  $E_v^1$, and $E_v^2$, where $E_v^l$ is the set of loops incident with $v$ which are distinct from $e$, $E_v^1 =E_H(v,N_v^1)$ (note that every edge in $E_v^1$ is a multi-edge.), and  $E_v^2 =E_H(v,N_v^2)$. Now we orient all edges in $E_{H_e}^-$. First, for every $s\in N_v^1$ we choose two edges belonging to $E_H(v,s)$, say $f^s$ and $g^s$. Let $A_E$ be the set of arcs obtained from $E_{H_e}^-$ by assigning any orientation to every loop in $E_v^l$, every edge in $E_v^2$ is oriented toward a vertex in $N_v^2$, while edges belonging to $E_v^1$ are oriented in such a~way that for every vertex $s\in N_v^1$ one chosen edge joining $v$ and $s$, say $f^s$, is oriented from $s$ to $v$, and all other edges belonging to $E_H(v,s)\setminus \{f^s\}$ are oriented toward $s$, see Fig. \ref{Rys-remark 1}. Let ${\cal P}_v$ be the family of oriented paths that consists of oriented 1-cycles $(v,h_A,v)$ (for every $h\in E_v^l$), oriented 2-cycles $(v,g_A^s, s,f_A^s,v)$ (for every $s\in N_v^1$), oriented 1-paths $(v,k_A,x)$ (for every $x\in N_v^2$ and every $k\in E_H(v,x)$), and $(v,l_A,y)$ (for every $y\in N_v^1$ and every $l\in E_H(v,y)\setminus \{f^y,g^y\}$). Finally, let $F_v$ be the digraph with vertex set $N_H[v]$ and arc set $A_e$. From the choice of ${\cal P}_v$ one can readily observe that $F_v$ and  ${\cal P}_v$ have the properties (1)--(4) stated in the definition of a good subgraph. Consequently, $H_e$ is a good subgraph in $H$.  \end{proof}

\begin{proposition}
\label{remark2} If $e$ is an edge joining two vertices, say $v$ and $u$, in a connected graph $H$,  then the subgraph $H_e$ of $H$ generated by $e$ is a good subgraph in $H$ if and only if  $H\not\in \{C_2, C_3\}$ and neither $v$ nor $u$ is adjacent to a pendant edge in $H$.
\end{proposition}
\begin{proof} It is easy to observe that if $H_e$ is a good subgraph in $H$, then $H\ne  K_2^2$, $H\ne  K_3$, and, by Observation \ref{leaves-supports}, neither $v$ nor $u$ is adjacent to a leaf in $H$.

Thus assume that  $H\ne  K_2^2$, $H\ne  K_3$, $e$ is an edge joining vertices $v$ and $u$ in $H$, and neither $v$ nor $u$ is adjacent to a leaf in $H$. We shall prove that $H_e$ is a good subgraph in $H$.  We consider two cases, namely $N_H(\{v,u\})= \{v,u\}$ and $\{v,u\} \varsubsetneq N_H(\{v,u\})$.

{\it Case 1. $N_H(\{v,u\})= \{v,u\}$.} In this case, $H$ is a graph of order~$2$. If $e$ is the only edge joining $v$ and $u$ in $H$, then $H_e$ is a~good subgraph in $H$ if and only if each of the vertices $v$ and $u$ is incident with a~loop in $H$. If $v$ and $u$ are joined by two parallel edges in $H$, then $H_e$ is a~good subgraph in $H$ if and only if at least one of the vertices $v$ and $u$ is incident with a~loop in $H$ (or, equivalently, $H_e$ is not  a good subgraph in $H$ if $H=K_2^2$).  Finally, if $v$ and $u$ are joined by at least three parallel edges in $H$, then $H_e$ is always a good subgraph in $H$. In every case it is straightforward to recognize arcs or directed paths forming the families of directed paths ${\cal P}_v$ and ${\cal P}_u$ and desired digraphs $F_v$ and $F_u$ in graphs $H_1,\ldots,H_6$ shown in Fig. \ref{Rys-remark 2}. (The digraphs $F_v$ and $F_u$ in Fig. \ref{Rys-remark 2} are drawn in blue and red, respectively.)

\begin{figure}[h!] \begin{center} \bigskip
{\epsfxsize=3.0in \epsffile{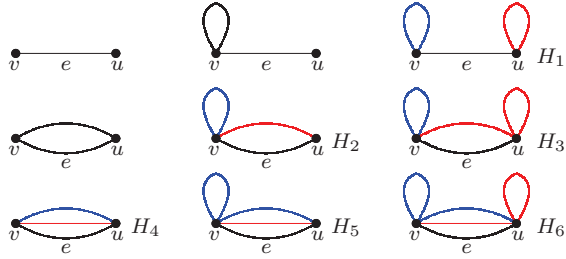}}

\vspace{-2mm}
\caption{Examples to Case 1} \label{Rys-remark 2}
\end{center}\end{figure}

\begin{figure}[t!] \begin{center} \bigskip
{\epsfxsize=4.0in \epsffile{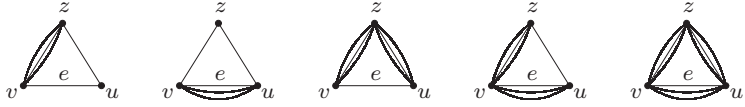}}

\vspace{-1mm}
\caption{An example to Observation \ref{remark2}} \label{Rys-remark 2koniec}
\end{center}\end{figure}

{\it Case 2. $\{v,u\} \varsubsetneq N_H(\{v,u\})$.} In this case, the set $N_{vu}^-=N_H(\{v,u\})\setminus \{v,u\}$ is nonempty, and, by our assumption, no vertex belonging to $N_{vu}^-$ is a leaf in $H$. The set $N_{vu}^-$ consists of five subsets:
$N_v^1= \{x\in N_{vu}^-\colon N_H(x)=\{v\}\}$,
$N_u^1= \{x\in N_{vu}^-\colon N_H(x)=\{u\}\}$,
$N_{vu}^1= \{x\in N_{vu}^-\colon N_H(x)=\{v,u\}\}$,
$N_v^2= \{x\in N_{vu}^-\colon \{v\} \varsubsetneq N_H(x)\}$, and
$N_u^2= \{x\in N_{vu}^-\colon \{u\} \varsubsetneq N_H(x)\}$.
The sets $N_v^1$, $N_u^1$, $N_{vu}^1$, and $N_v^2 \cup N_u^2$ are disjoint, and it is possible that some of them are empty. Let $A_E$ be a set of arcs obtained by assigning an orientation to every edge belonging to the set $E_{H_e}^-$, that is, to every edge incident with $v$ or $u$ and different from $e$. The set $A_E$ and families ${\cal P}_v$ and ${\cal P}_u$ of directed paths that begin at $v$ and $u$, respectively, are defined in the following way: \\ [-24pt]
\begin{enumerate}
\item[$(1)$] To every loop $f$ incident with $x$ we assign an arbitrary orientation $f_A$, and add the 1-cycle $(x,f_A,x)$ to ${\cal P}_x$ if $x\in \{v,u\}$. \1
\item[$(2)$] If an edge $f$ belongs to $E_H(\{v,u\},N_v^2\cup N_u^2)$, then by $f_A$ we denote the orientation of $f$ toward $N_v^2\cup N_u^2$. In addition, if $\varphi_H(f)=\{x,y\}$, where $x\in \{v,u\}$ and $y\in N_x^2$, then the 1-path $(x,f_A,y)$ we add to ${\cal P}_x$. \1
\item[$(3)$] If $x\in N_v^1\cup N_u^1$, and $y$ is the only neighbor of $x$ (belonging to $\{v,u\}$), then (as in the proof of Proposition \ref{remark1}) we choose two edges belonging to $E_H(x,y)$, say $f^x$ and $g^x$, one of them, say $f^x$, obtain an orientation from $x$ to $y$, and all other edges belonging to $E_H(x,y)\setminus \{f^x\}$ are oriented toward $x$. In this case, the 2-cycle $(y,g_A^x,x,f_A^x,y)$ and all the 1-paths $(y,h_A,x)$, for every $h\in E_H(x,y)\setminus \{f^x,g^x\}$ (if this set is nonempty), we add to ${\cal P}_y$. \1
\item[$(4)$] If the set $N_{vu}^1$ is nonempty, then we distinguish two cases:
\begin{enumerate}
\item If at least one of the sets $E_H(v,N_v^1\cup N_u^1)\cup E_v^l$ and $E_H(u,N_v^1\cup N_u^1)\cup E_u^l$ is nonempty, say $E_H(v,N_v^1\cup N_u^1)\cup E_v^l \not =\emptyset$, then for every $z\in N_{vu}^1$, we choose two edges belonging to $E_H(z,\{v,u\})$, say $f^z\in E_H(z,v)$ and $g^z \in E_H(z,u)$, orient $f^z$ toward $v$, all other edges belonging to $E_H(z,\{v,u\})\setminus \{f^z\}$ are oriented toward $z$, and to every edge in $E_H(v,u)\setminus \{e\}$ (if this set is nonempty) we choose an arbitrary orientation, say from $u$ to $v$. Now the $2$-path $(u,g_A^z,z,f_A^z,v)$, the $1$-paths $(u,h_A,z)$ (for every $h\in E_H(u,z)\setminus \{g^z\}$) and $(u,k_A,v)$ (for every $k\in E_H(u,v)\setminus \{e\}$) we add to ${\cal P}_u$, while the 1-paths
    $(v,l_A,z)$ (for every $l\in E_H(v,z)\setminus \{f^z\}$) to ${\cal P}_v$. \1
\item If both the sets $E_H(v,N_v^1\cup N_u^1)\cup E_v^l$ and $E_H(u,N_v^1\cup N_u^1)\cup E_u^l$ are empty, then we consider two subcases: \begin{enumerate}
    \item[(b1)] If $|N_{vu}^1|\ge 2$, and if $C$ is a smallest subset of $E_H(\{v,u\},N_{vu}^1)$ that covers the vertices in $\{v,u\}\cup N_{vu}^1$, then we orient the edges in $C$ toward $\{v,u\}$, the edges in $E_H(\{v,u\},N_{vu}^1)\setminus C$ toward $N_{vu}^1$, and the edges belonging to $E_H(v,u)\setminus \{e\}$ (if this set is nonempty) in an arbitrary way, again say from $u$ to $v$. We may assume that $N_{vu}^1=\{z_1,\ldots,z_k\}$, $C=\{f^{z_1},\ldots, f^{z_k}\}$, and $\varphi_H(f^{z_i})= \{z_i,x_i\}$, where $x_i\in \{v,u\}$ for $i\in [k]$. Let $D=\{g^{z_1}, \ldots, g^{z_k}\}$, where  $\varphi_H(g^{z_i})= \{z_i,y_i\}$ and $y_i$ is the only element of $\{v,u\}\setminus \{x_i\}$ for $i\in [k]$. Now, all 1-paths $(u,l_A,v)$  (for every $l\in E_H(v,u)\setminus \{e\}$) and 1-paths $(u,l_A,z_i)$ (if $i\in [k]$ and  $l\in E_H(u,z_i) \setminus (C\cup D)$) we add to ${\cal P}_u$, while 1-paths $(v,p_A,z_i)$ (if $i\in [k]$ and  $p\in E_H(v,z_i) \setminus (C\cup D)$) we add to ${\cal P}_v$. Finally, the 2-path $(y_i,g_A^{z_i},z_i, f_A^{z_i},x_i)$, $i\in [k]$, we add to ${\cal P}_v$ (${\cal P}_u$, resp.) if and only if $y_i=v$ ($y_i=u$, resp.).\1
    \item[(b2)] If $|N_{vu}^1|=1$, say $N_{vu}^1= \{z\}$, then, since $V_H=\{v,u,z\}$ and $H\ne  K_3$, $H$ is a proper spanning supergraph of $K_3$ and, therefore, it has parallel edges (as $E_v^l = E_u^l= \emptyset$). Without loss of generality, we assume that $v$ is incident with parallel edges. There are five cases to consider, and they are sketched in Fig. \ref{Rys-remark 2koniec}. In each of these cases, let $f^z$ and $g^z$ be an edge belonging to $E_H(v,z)$ and $E_H(u,z)$, respectively. We orient $f^z$ toward $v$, all other edges belonging to $E_H(\{v,u\},z)\setminus \{f^z\}$ we orient toward $z$, and the edges belonging to $E_H(v, u)\setminus \{e\}$ (if $E_H(v, u) \setminus \{e\} \ne  \emptyset$) are directed toward $u$. Now, the 2-path $(u,g_A^z,z,f_A^z,v)$ and 1-paths $(u,h_A,z)$ ($h\in E_H(u,z)\setminus \{g^z\}$) form the family ${\cal P}_u$, while 1-paths $(v,l_A,z)$ ($l\in E_H(v,z)\setminus \{f^z\}$) and $(v,p_A,u)$ ($p\in E_H(v,u)\setminus \{e\}$) form the family ${\cal P}_u$.
\end{enumerate}\end{enumerate}\end{enumerate}

Let $F_v$ and $F_u$ be digraphs generated by arcs belonging to families ${\cal P}_v$ and ${\cal P}_u$, respectively. Since families ${\cal P}_v$ and ${\cal P}_u$ consist of 1- and 2-paths, we observe that the digraphs $F_v$ and $F_u$, and families ${\cal P}_v$ and ${\cal P}_u$ have the properties (1)--(4) stated in the definition of a good subgraph. Consequently, $H_e$ is a good subgraph in $H$.
 \end{proof}

\begin{remark}  Let $H$ be a graph without isolated vertices and let $I$ be a proper subset of $V_H$. Then the induced subgraph $H[I]$ is a good subgraph in $H$ if $1\le d_{H[I]}(v)< d_H(v)$ for every $v\in I$, and $N_H(x)\setminus I\ne \emptyset$ for every $x\in N_H(I)\setminus I$. {\rm A proof of this statement is similar to the proofs of Propositions~\ref{remark1} and~\ref{remark2} and is omitted.} \end{remark}

As a consequence of Observation \ref{leaves-supports} and Propositions~\ref{remark1} and~\ref{remark2}, we have the following two corollaries.

\begin{cor2} \label{corollary-small-good-subgraphs} \label{corollary-2-small-good-subgraphs} A connected graph has a good subgraph if and only if it has a good subgraph generated by a loop or by an edge, that is, if and only if $H\not\in\{C_1, C_2,C_3\}$ and it has an edge $($or a~loop$)$ which is neither a pendant edge nor adjacent to a~pendant edge in $H$.  \end{cor2}

\begin{cor2} \label{good subgraphs in trees} A tree $T$ of order at least~$2$ has a good subgraph if and only if it has an edge which is neither a pendant edge nor adjacent to a pendant edge in $T$. \end{cor2}

It might be thought that a $2$-subdivision graph of a graph having a good subgraph is not a~minimal DTDP-graph, but it is not true in general since, for example, $K_1^2$ has a good subgraph and its $2$-subdivision graph $S_2(K_1^2,{\cal P}_4)$ shown in Fig. \ref{examples-K-1-m} is a~mini\-mal DTDP-graph. For this reason in the next theorem (which is important in our characterization of the minimal DTDP-graphs) we only consider $2$-subdivision graphs without loops, that is, $2$-subdivision graphs in which no twin parts corresponding to a loop are contracted into a single edge and, in consequence, forming a loop in the $2$-subdivision graph.

\begin{thm} 
\label{zakazane-w-minimalnych-DTDP} Let $H$ be a connected graph without isolated vertices, and let ${\cal P}=\{{\cal P}(v)\colon v\in V_H\}$ be a family in which ${\cal P}(v)$ is a~partition of the neighborhood $N_{S_2(H)}(v)$ for every  $v\in V_H$. If $H$ has a good subgraph and ${\cal P}$ does not contract in $S_2(H,{\cal P})$ any twin parts corresponding to a loop in $H$, then the $2$-subdivision graph $S_2(H,{\cal P},\theta)$ is a non-minimal DTDP-graph $($for every positive function $\theta \colon L_{S_2(H,{\cal P})}\to \mathbb{N}$$)$. \end{thm}
\begin{proof} Let $Q$ be a good subgraph in $H$. By Corollary \ref{corollary-small-good-subgraphs} we may assume that $Q$ is a good subgraph generated by a loop or by an edge. By Proposition~\ref{observ-S2(H)-jest-DPDP-grafem} the $2$-subdivision graph $S_2(H,{\cal P},\theta)$ is a~DTDP-graph. We shall prove that $S_2(H,{\cal P},\theta)$ is not a minimal DTDP-graph. By Observation \ref{observ-S2(H)-stron-weak-supports} it suffices to show that  $S_2(H,{\cal P})$ is a non-minimal DTDP-graph. By virtue of Proposition~\ref{o sklejaniu krawedzi} it suffices to show this non-minimality only in the case when ${\cal P}$ contracts in $S_2(H,{\cal P})$ far parts of adjacent pendant edges in $H$ (since we have assumed that
${\cal P}$ does not contract in $S_2(H,{\cal P})$ any twin parts corresponding to any loop in $H$). In such case it is possible to observe that $S_2(H,{\cal P})$ is a non-minimal DTDP-graph if and only if $S_2(H)$ is a~non-minimal DTDP-graph. Thus it remains to prove that $S_2(H)$ is a non-minimal DTDP-graph.

Assume first that the good subgraph $Q$ in $H$ is generated by a loop, say by a~loop $e$ incident with a vertex $v$. It is obvious that $H=K_1^s$ has a good subgraph and $S_2(K_1^s)$ is a non-minimal DTDP-graph if and only if $s\ge 2$. Thus assume that $H$ is a~connected graph of order at least~$2$. For simplicity, as far as possible, we adopt the notation from the proof of Proposition \ref{remark1}.
For ease of presentation, we assume that $N_v^1= \{v^1,\ldots,v^k\}$ and
$E_H(v,v^i)= \{e_i^1,\ldots, e_i^{j_i} \}$ (where $j_i\ge 2$ as every edge in $E_H(v,N_v^1)$ is a multi-edge) for every $v^i\in N_v^1$. We may assume that $A_E$ is an orientation of $E_{Q}^-$ (of the set of edges or loops belonging to $E_H\setminus E_{Q}= E_H\setminus \{e\}$ that are incident with $v$) such that every loop belonging to $E_v^l$ obtain an arbitrary direction, every edge belonging to $E_v^2$ is directed toward $N_v^2$, and edges belonging to $E_v^1$ are oriented in such a way that for every vertex $v^i\in N_v^1$ the edge $e_i^1$ is directed from $v^i$ to $v$, and all other edges belonging to $E_H(v,v^i)$ are directed toward $v^i$. Let $F_v$ be the digraph generated by the arcs belonging to $A_E$. Let ${\cal P}_v$ be the family consisting of all directed 2-cycles $(v,e_i^2,v^i,e_i^1,v)$ (for $i \in [k]$) and of all directed 1-paths (and 1-cycles) generated by all other arcs of $F_v$, see the left part of Fig. \ref{Rys-dowod-thm 4.5}. The digraph $F_v$ and the family ${\cal P}_v$, as in the proof of Proposition \ref{remark1}, have the properties (1)--(4) stated in the definition of a good subgraph, implying that $Q$ is a desired good subgraph in $H$.

Let $G'$ be the proper spanning subgraph obtained from $S_2(H)$ by removing the ``middle" edge $v_e^1v_e^2$ from the $3$-cycle corresponding to the loop $e$ of $Q$, and the third edge from the $4$-path corresponding to the last arc in every directed path in ${\cal P}_v$, as illustrated in the right part of Fig. \ref{Rys-dowod-thm 4.5}. Formally, $G'$ is the proper spanning subgraph  of $S_2(H)$ with edge set
\[
E_{G'}= E_{S_2(H)}\setminus \left( \{v_e^1v_e^2\} \cup \{vv_f^2\colon f\in E_v^l\} \cup \bigcup_{u\in N_v^2}\{ uu_g\colon g\in E_H(v,u)\} \cup \bigcup_{i=1}^k \{vv_{e_i^1},
v^iv_{e_i^3}^i,\ldots, v^iv_{e_i^{j_i}}^i\}\right).
\]

All that remains to prove is that $G'$ is a DTDP-graph. It suffices to observe that every component of $G'$ is a $2$-subdivision graph.
Let $G_v'$ be the component of $G'$ containing the vertex $v$. We note that $G_v'$ belongs to the family $\cal F$ and, therefore, it is a DTDP-graph by Corollary \ref{Corollary family F}.
If the set $N_v^2$ is empty, then $G'=G_v'$ and the desired result follows. Thus assume that the set $N_v^2$ is non-empty. Then, since every edge belonging to the set
\[
\bigcup_{u\in N_v^2}\{ uu_g \colon g\in E_H(v,u)\}
\]
joins a vertex in $N_v^2$ to a vertex in $V_{G_v'}$, while every edge belonging to the set
\[
\{v_e^1v_e^2\} \cup \{vv_f^2\colon f\in E_v^l\} \cup \bigcup_{i=1}^k \{vv_{e_i^1},v^iv_{e_i^3}^i,\ldots, v^iv_{e_i^{j_i}}^i\}
\]
joins two vertices belonging to $V_{G_v'}$, the subgraph $G''=G'-V_{G_v'}$ is an induced subgraph of $G'$ and, in addition, $G''$ is a $2$-subdivision graph, $G''=S_2(H-V_{G_v'})$. Thus, by Proposition~\ref{observ-S2(H)-jest-DPDP-grafem}, $G''$ is DTDP-graph. Consequently, since $G_v'$ and $G''$ are DTDP-graphs, the proper spanning subgraph $G'$ of $S_2(H)$ is a DTDP-graph and, therefore,  $S_2(H)$ is a non-minimal DTDP-graph.
(For example, in Fig. \ref{Rys-dowod-thm 4.5} it is $G'=S_2(C_4)\cup S_2(K_{1,10}^2,{\cal P},\theta)$, where $K_{1,10}^2$ is a~star $K_{1,10}$ with two subdivided edges, ${\cal P}$ in $S_2(K_{1,10}^2)$ contracts all ten neighbors of the vertex corresponding to the central vertex of $K_{1,10}$ (or $K_{1,10}^2$), and finally the ``new" pendant edge in $S_2(H,{\cal P})$ is replaced by twin pendant edges (using the function $\theta$).))

\begin{figure}[h!] \begin{center} \bigskip
{\epsfxsize=5.35in \epsffile{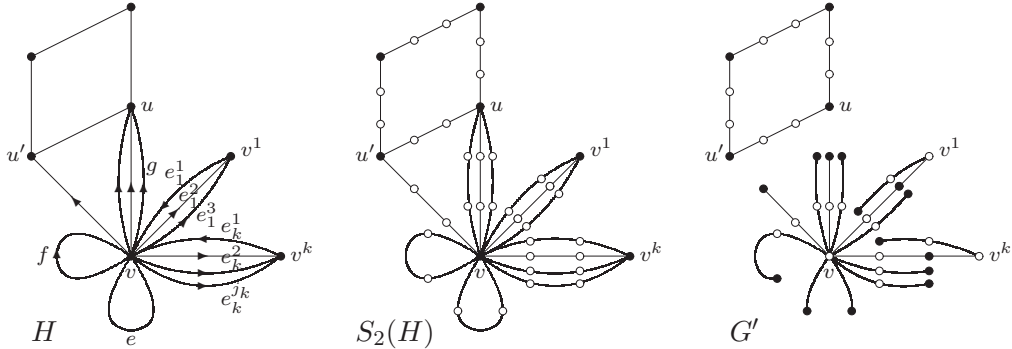}}

\vspace{-2mm}
\caption{Graphs $H$, $S_2(H)$, and a spanning subgraph $G'$ of $S_2(H)$} \label{Rys-dowod-thm 4.5}
\end{center}\end{figure}

Assume now that the good subgraph $Q$ in $H$ is generated by an edge, say by an edge $e$ which joins two vertices $v$ and $u$ in $H$. We know that $S_2(H)$ is a~DTDP-graph and we shall prove that $S_2(H)$ is a non-minimal DTDP-graph. We already know that $S_2(H)$ is a non-minimal DTDP-graph if $H$ has a good subgraph generated by a~loop. Thus assume that no subgraph of $H$ generated by a loop is a good subgraph in $H$. Consequently, since $H$ is a connected graph of order at least 2, every loop in $H$ is incident with a support vertex, and, in particular, neither $v$ nor $u$ is incident with a loop. Certainly, neither $v$ nor $u$ is a support vertex in $H$. As in the proof of Proposition \ref{remark2} we consider two cases, namely $N_H(\{v,u\})= \{v,u\}$ and $\{v,u\} \varsubsetneq N_H(\{v,u\})$.

{\it Case 1. $N_H(\{v,u\})= \{v,u\}$.} In this case, since neither $v$ nor $u$ is incident with a loop,  $H=K_2^s$ and $s\ge 1$. From the fact that $K_2^s$ has a good subgraph it follows that $s\ge 3$. Certainly,  $S_2(K_2^s)$ is a non-minimal DTDP-graph if $s\ge 3$.

\begin{figure}[h!] \begin{center} \bigskip
{\epsfxsize=5.75in \epsffile{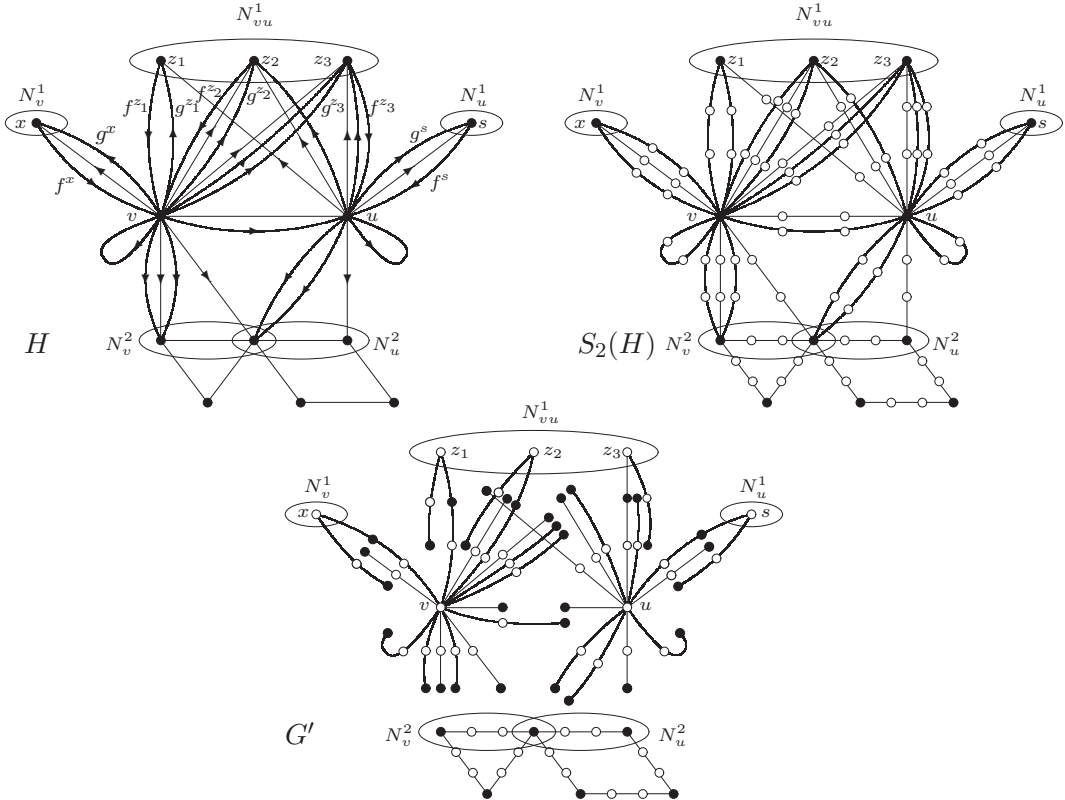}}

\vspace{-2mm}
\caption{Graphs $H$, $S_2(H)$, and a spanning subgraph $G'$ of $S_2(H)$ (Case 2)} \label{Rys-dowod-thm 4.5-case 2}
\end{center}\end{figure}

{\it Case 2. $\{v,u\} \varsubsetneq N_H(\{v,u\})$.} For simplicity we use the same notation as in the second part of the proof of Proposition \ref{remark2}. Consider the orientation $A_E$ of $E_Q^-$, the families ${\cal P}_v$ and ${\cal P}_u$, and the digraphs $F_v$ and $F_u$, introduced in Case 2 of the proof of Proposition \ref{remark2}. Let $G'$ be the spanning subgraph of $S_2(H)$  obtained from $S_2(H)$ by removing the middle edge $v_eu_e$ from the 4-path $(v,v_e,u_e,u)$ corresponding to the edge $e$, and the third edge from each 4-path corresponding to the last arc in every directed path in ${\cal P}_v$ or ${\cal P}_u$, see the lower part of Fig. \ref{Rys-dowod-thm 4.5-case 2}. As in the first part of the proof, $G'$ is a DTDP-graph. Consequently, the proper spanning subgraph $G'$ of $S_2(H)$ is a~DTDP-graph and therefore $S_2(H)$ is a non-minimal DTDP-graph.
\end{proof}

It follows from Corollary \ref{corollary-2-small-good-subgraphs} that every path $P_n$ (with $n\ge 6$) and every cycle $C_m$ (with $m\ge 4$) has a~good subgraph, and, therefore, Proposition~\ref{observ-S2(H)-jest-DPDP-grafem}, Theorem \ref{zakazane-w-minimalnych-DTDP}, and a simple verification justify the following remark about minimal  DTDP-paths and minimal DTDP-cycles, see Observation \ref{examples-plus}.

\begin{remark} \label{rem-mini-cycle-path}  If $C_m$ is a cycle of size $m$, then $S_2(C_m)$ is a DTDP-graph for every positive integer $m$, but $S_2(C_m)$ is a minimal DTDP-graph if and only if $m\in \{1, 2, 3\}$. If $P_n$ is a path of order $n$, then $S_2(P_n)$ is a DTDP-graph for every integer $n\ge 2$, while $S_2(P_n)$ is a minimal DTDP-graph if and only if $n\in \{2, 3, 4, 5\}$.  \end{remark}

\section{Structural characterization of the DTDP-graphs}%

The next theorem presents general properties of \DT-pairs in  minimal DTDP-graphs.

\begin{thm} 
\label{observ-2-minimal-DTDP-graphs} Let $G$ be a~connected minimal DTDP-graph, and let $(D,T)$ be a~\DT-pair in $G$. Then the following four statements hold: \\ [-24pt]
\begin{enumerate}
\item[$(1)$] The set $D$ is a~maximal independent set in $G$.
\item[$(2)$] Every component of $G[T]$ is a star or it is a graph of order~$1$ and size~$1$.
\item[$(3)$] If $x\in T$, then  $|N_G(x)\setminus T|=1$ or $N_G(x)\setminus T$ is a~nonempty subset of $L_G$. In addition, if $x$ is a leaf in a star of order at least 3 in $G[T]$, then $N_G(x)\setminus T$ is a~nonempty subset of $L_G$.
\item[$(4)$] \label{propertis-of-minimal} The graph $G$ is a $2$-subdivision graph $S_2(H,{\cal P},\theta)$  of some connected graph $H$, where ${\cal P}=\{{\cal P}(v)\colon v\in V_H\}$ is a family in which ${\cal P}(v)$ is a partition of the neighborhood $N_{S_2(H)}(v)$ for $v\in V_H$, and $\theta \colon L_{S_2(H,{\cal P})}\to \mathbb{N}$ is a positive function.
\end{enumerate}
\end{thm}
\begin{proof} (1) If the set $D$ were not independent, then two vertices belonging to $D$, say $x$ and $y$, would be adjacent, and then $(D,T)$ would be a~\DT-pair in $G-xy$,  contradicting the minimality of $G$. Now, since $D$ is both an independent set and a dominating set of $G$, $D$ is a~maximal independent set in $G$.

(2) Since $T$ is a total dominating set of $G$, by definition, $G[T]$ has no isolated vertex. Consequently, every component of $G[T]$ is either of order~$1$ and size at least~$1$ or has order at least~$2$. From the minimality of $G$, a~component of order~$1$ in $G[T]$ has exactly one loop incident with its only vertex. Now let $F$ be a component of order at least~$2$ in $G[T]$. To prove that $F$ is a star, it suffices to show that if distinct vertices are adjacent in $G[T]$, then at least one of them is a leaf in $G[T]$. If  $x$ and $y$ are adjacent in $G[T]$ and neither of them is a~leaf in $G[T]$, then $(D,T)$ would be a~\DT-pair in $G-xy$, violating the minimality of $G$.

(3) Assume that $x\in T$. Then $N_G(x)\setminus T$ is a nonempty subset of $D$ (since $D=V_G\setminus T$ is a~dominating set in $G$). Therefore, since $L_G\subseteq D$ (by Observation \ref{observ-1}), either $N_G(x)\setminus T$ is a nonempty subset of $L_G$ or $N_G(x)\setminus (L_G\cup T)$ is nonempty. It remains to prove that if $N_G(x)\setminus (L_G\cup T)$ is nonempty, then $|N_G(x)\setminus T|=1$. Assume that $y\in N_G(x)\setminus (L_G\cup T)$. Then, since $D$ is independent and $y\in D\setminus L_G$, the set $N_G(y)\setminus \{x\}$ is a nonempty subset of $T$, say $x'\in N_G(y)\setminus \{x\}$.  Now suppose that $|N_G(x)\setminus T|\ge 2$, and let $y'$ be any vertex in $(N_G(x)\setminus T)\setminus \{y\}$. Then, since $x$ is dominated by $y'$ and $y$ is totally dominated by $x'$,  the pair $(D,T)$ is a~\DT-pair in $G-xy$, contradicting the minimality of $G$.

Assume now that $x$ is a leaf in a star of order at least~$3$ in $G[T]$. Let $x'$ be the only neighbor of $x$ in $G[T]$, and let $y\in N_G(x)\setminus \{x'\}$. It remains to show that $y$ is a leaf in $G$. Suppose that $y$ is not a leaf in $G$. Then $N_G(y)\setminus \{x\}\ne  \emptyset$, and, if $x''\in N_G(y)\setminus \{x\}$, then the pair $(D\cup\{x\}, T\setminus \{x\})$ is a~\DT-pair in $G-xy$, a~contradiction which completes the proof of (3).

(4) Let $G=(V_G, E_G, \varphi_G)$ be a graph (where, as usually, $\varphi_G(e)$ denotes the set of vertices incident with $e\in E_G$). Assume that $G$ is a~minimal DTDP-graph, and let $(D,T)$ be a~\DT-pair in $G$. The minimality of $G$ implies that $G$ has neither multi-edges nor multi-loops. With respect to Observation \ref{observ-S2(H)-stron-weak-supports} we may assume that $G$ has no strong supports. This assumption and the already proved property (3) imply that every vertex $v$ belonging to $T$ has exactly one neighbor in $D$, and, in addition, this unique neighbor of $v$ is a leaf in $G$ if $v$ is a leaf of a star of order at least 3 in $G[T]$. Consequently, if $e$ is an edge in $G[T]$, then the subset $N_G(\varphi_G(e))\setminus T$ of $D$ is of order 1 or 2. This implies that the triple $H=(V_H, E_H,\varphi_H)$ in which $V_H=D$, $E_H= E_{G[T]}$, and $\varphi_H \colon E_H\to 2^{V_H}$ is a function such that $\varphi_H(e)= N_{G}(\varphi_G(e))\setminus T$ for each edge $e\in E_H$, is a~well-defined graph with possible multi-edges or multi-loops. (If $e$ is an edge in $G[T]$ and $\varphi_H(e)=\{a,b\}$, then $e$ is an edge which joins the vertices $a$ and $b$ in $H$. Similarly, if $e$ is an edge or a loop in $G[T]$ and $\varphi_H(e)=\{a\}$, then $e$ is a~loop which joins $a$ to itself in $H$.)
Now, to restore the graph $G$ from the multi-graph $H$, we first form $S_2(H)$ inserting two new vertices into each edge and each loop of $H$. More precisely, if an edge $e$ joins vertices $a$ and $b$ in the multi-graph $H$ (that is, if $\varphi_H(e)=\{a,b\}$), then by $a_e$ and $b_e$ we denote the two mutually adjacent vertices inserted into the edge $e$, where $a_e$ and $b_e$ are adjacent to $a$ and $b$, respectively. (If $e$ is a loop incident with a vertex $a$ (that is, if $\varphi_H(e)=\{a\}$), then by $a_{e}^1$ and $a_{e}^2$ we denote two mutually adjacent vertices inserted into the loop $e$ and both adjacent to $a$.) Let ${\cal P}=\{{\cal P}(v)\colon v\in V_H\}$ be a family in which the partition ${\cal P}(v)$ of the set $N_{S_2(H)}(v)$ (for $v\in V_H$ ($\subseteq V_{S_2(H)}$) is defined as follows:

\begin{itemize}
\item If $e$ is a loop in $G[T]$ incident with a vertex $N_G(v)$, then we let the 2-element set $\{v_e^1, v_e^2\}$ be an element of ${\cal P}(v)$.
\item If $e$ is an edge (not a loop) in $G[T]$ and $\varphi_G(e)\subseteq N_G(v)$, then we choose both one-element sets $\{v_e^1\}$  and $\{v_e^2\}$ to belong to ${\cal P}(v)$.
\item If $\{e_1,\ldots, e_k\}$ is the edge set of a star in $G[T]$ and the central vertex of this star is in $N_G(v)$ (or exactly one vertex  of this star is in $N_G(v)$, if $k=1$), then we select the set $\{v_{e_1},\ldots, v_{e_k}\}$ as an element of ${\cal P}(v)$.
\end{itemize}

From the above definition of the family ${\cal P}$, the $2$-subdivision graph $S_2(H,{\cal P})$ ($=S_2(H,{\cal P},\theta)$ if $\theta(x)=1$ for every $x\in L_{S_2(H,{\cal P})}$) obtained from $S_2(H)$ is isomorphic to the graph $G$, see Fig. \ref{grafy G, H i S_2(H)} for an illustration. \end{proof}

\begin{figure}[h!] \begin{center} \bigskip
{\epsfxsize=6.5in \epsffile{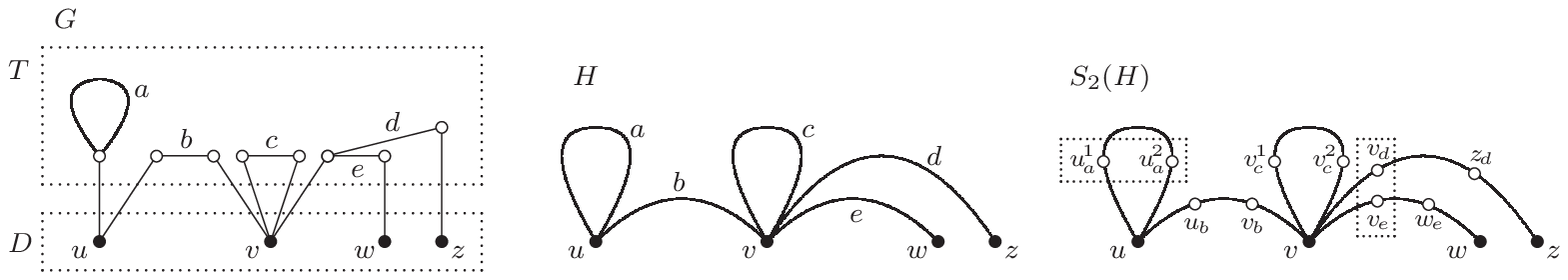}}

\vspace{-2mm}
\caption{} \label{grafy G, H i S_2(H)}
\end{center}\end{figure}

From Theorem \ref{observ-2-minimal-DTDP-graphs}\,(4) every minimal DTDP-graph is a $2$-sub\-di\-vi\-sion graph of some graph. The converse, however, is not true in general. It is easy to check that neither of the $2$-subdivision graphs $S_2(H)$, $S_2(H,{\cal P})$, and $S_2(H,{\cal P},\theta)$ presented  in Fig. \ref{2-subdivision-graph-new} is a minimal DTDP-graph. In our last theorem we present the main structural characterization of minimal DTDP-graphs without loops.

\begin{thm}
\label{minimal-DTDP-graphs-2} Let $G$ be a connected graph of order at least 3 that has no loops. Then the following statements are equivalent:\\ [-24pt]
\begin{enumerate}
\item[$(1)$] The graph $G$ is a minimal DTDP-graph; \2
\item[$(2)$] Either $\rm (2a)$ $G\in \{C_3, C_6, C_9\}$ or $\rm (2b)$ $G$ is a $2$-subdivision graph, say $G=S_2(H,{\cal P},\theta)$ $($where $H$ is a connected graph of order at least~$2$, ${\cal P}=\{{\cal P}(v)\colon v\in V_H\}$ is a~family in which ${\cal P}(v)$ is a partition of the set  $N_{S_2(H)}(v)$ for $v\in V_H$, and $\theta \colon L_{S_2(H,{\cal P})}\to \mathbb{N}$ is a positive function$)$ in which $\rm (2b1)$ the pair $(V_G^o, V_G^n) =(V_{S_2(H,{\cal P},\theta)}^o,V_{S_2(H,{\cal P},\theta)}^n)$ is the only \DT-pair and, in addition, in which $\rm (2b2)$ every component of $G[V_{G}^n]$ is a star. \2
\item[$(3)$] The graph $G$ is a $2$-subdivision graph, say $G=S_2(H,{\cal P},\theta)$, where either
$\rm (3a)$ $H\in \{C_1, C_2, C_3\}$ or $\rm (3b)$ $H$ is a connected graph of order at least 2 in which $\rm (3b1)$ every non-pendant edge $($and every loop$)$ is adjacent to a~pendant edge, $\rm (3b2)$ ${\cal P}=\{{\cal P}(v)\colon v\in V_H\}$ is a family in which every ${\cal P}(v)$ is a partition of the set  $N_{S_2(H)}(v)$ for every $v\in V_H$, and ${\cal P}$ contracts in $S_2(H,{\cal P})$ only far parts of adjacent pendant edges of $H$ $($if any$)$, and $\rm (3b3)$ $\theta \colon L_{S_2(H,{\cal P})} \to \mathbb{N}$ is a positive function.
\end{enumerate} \end{thm}
\begin{proof} $(1) \Rightarrow (2)$. Assume first that $G$ is a connected minimal DTDP-graph. Then $G$ has no multi-edges and Theorem \ref{observ-2-minimal-DTDP-graphs} implies that $G$ is a $2$-subdivision graph, i.e., $G=S_2(H,{\cal P},\theta)$ $($for some connected graph $H$ without a good subgraph (by Theorem \ref{zakazane-w-minimalnych-DTDP}), some family ${\cal P}=\{{\cal P}(v)\colon v\in V_H\}$ in which ${\cal P}(v)$ is a partition of the neighborhood $N_{S_2(H)}(v)$ $($for $v\in V_H$$)$ which contracts at most far parts of pendant edges (by Proposition~\ref{o sklejaniu krawedzi}), and a positive function $\theta \colon L_{S_2(H,{\cal P})}\to \mathbb{N}$$)$. By Observation \ref{observ-S2(H)-stron-weak-supports} we may  assume that $G$ has no strong supports, and therefore we may assume that $G=S_2(H,{\cal P})$ (as $S_2(H,{\cal P})$ and $S_2(H,{\cal P},\theta)$ are isomorphic if $\theta(x)=1$ for every $x\in L_{S_2(H,{\cal P})}$). By Proposition~\ref{observ-S2(H)-jest-DPDP-grafem}, the pair $(D,T)=(V_G^o, V_G^n)$ is a~\DT-pair in $G$. In addition, the minimality of $G$, Theorem \ref{observ-2-minimal-DTDP-graphs} and our assumption that $G$ has no loop imply that every component of $G[V_G^n]$ is a~star. Thus, it remains to prove that either $G$ is a cycle of length 3, 6 or 9, or the pair $(V_G^o, V_G^n)$ is the only \DT-pair in $G$.  We consider three cases depending on $\Delta(H)$.

{\it Case 1. $\Delta(H)=1$.} In this case, $H=P_2$ and $G=S_2(H,{\cal P})=P_4$ (as by our assumption $G$ has no strong supports). Moreover, $(V_G^o,V_G^n)= (L_G,S_G)$ is the only \DT-pair in $G$.

{\it Case 2. $\Delta(H)=2$.} In this case, either $H=C_m$ ($m\ge 1$) or $H=P_n$ ($n\ge 3$). But, since $S_2(H,{\cal P})$ is a minimal DTDP-graph, Remark~\ref{rem-mini-cycle-path} implies that either $H=C_m$ and $m\in \{1, 2, 3\}$, or $H=P_n$ and $n\in\{3, 4, 5\}$. Now, depending on ${\cal P}$, $S_2(C_1,{\cal P})=C_3$ or $S_2(C_1,{\cal P})=G_1$ (where $G_1$, shown in Fig. \ref{minimal-graphs-with-loops}, is the smallest DTDP-graph) and only $C_3$ has the desired properties (as $G_1$ has a loop). It is also a simple matter to observe that $S_2(P_3, {\cal P})=S_2(P_3)=P_7$ or $S_2(P_3,{\cal P})=P_3\circ K_1$ and each of these graphs has the desired properties. Simple verifications and Proposition~\ref{o sklejaniu krawedzi} show that each of the graphs $S_2(C_2,{\cal P})$ (see Fig. \ref{examples}), $S_2(C_3,{\cal P})$, $S_2(P_4,{\cal P})$, and $S_2(P_5,{\cal P})$ is a minimal DTDP-graph if and only if $S_2(C_2,{\cal P})=S_2(C_2)=C_6$,  $S_2(C_3,{\cal P})=S_2(C_3)=C_9$, $S_2(P_4,{\cal P})=S_2(P_4)=P_{10}$, and $S_2(P_5,{\cal P})=S_2(P_5)=P_{13}$, respectively. Certainly, each of these four graphs has the desired properties.

{\it Case 3. $\Delta(H)\ge 3$.} In this case we claim that $(D,T)= (V_G^o, V_G^n)$ is the only \DT-pair in $G$. Suppose, to the contrary, that $(D',T')$ is another \DT-pair in $G$. Then, since $D$ and $D'$ are maximal independent sets in $G$ (by Theorem \ref{observ-2-minimal-DTDP-graphs}) and $D\ne D'$, each of the sets $D\setminus D'$ and $D'\setminus D$ is a nonempty subset of $T'$ and $T$, respectively. Let $v$ be a vertex of maximum degree among all vertices in $D\setminus D' \subseteq T'$. Since $v \in T'$, it follows from Theorem \ref{observ-2-minimal-DTDP-graphs} that $d_H(v)\ge 2$. If $H=K_1^s$ and $s\ge 2$ (as $d_H(v)\ge 3$), then $K_1^1$ would be a good subgraph in $H$, which is impossible. Hence, $H$ has order at least~$2$. We consider two possible cases.

{\it Case 3.1. There is a loop, say $e$, at $v$.} In this case, a pendant edge, say $f$, is incident with $v$, as otherwise, by Proposition \ref{remark1}, the subgraph generated by $e$ would be a good subgraph in $H$, which again is impossible.  Assume that $f$ joins the vertex $v$ with a leaf $u$ in $H$. We claim that $v$ belongs to the set $D'$. Let $A$ be the only set in ${\cal P}(v)$ which contains the vertex $v_f$. By Observation \ref{observ-1}, $u\in D'$ and $(u_f,\{u_f\}) \in T'$. Thus, $(v,A)\in T'$ as $T'$ is a total dominating set of $G$ and $(v,A)$ is the only neighbor of $(u_f,\{u_f\})$ in $T'$. Finally, this implies that $v\in D'$ as $D'$ is a~dominating set of $G$ and $v$ is the only neighbor of $(v,A)\in T'=V_G\setminus D'$. Consequently, $v\in D'$ and  $v\in D\setminus D'$ (by the choice of $v$), a contradiction.

{\it Case 3.2. No loop is incident with $v$.} In this case, let $f_1, \ldots, f_k$ be the edges incident with $v$ in $H$, say $\varphi_H(f_i)=\{v,v^i\}$ for $i \in [k]$ ($k\ge 2$). If at least one of the edges $f_1, \ldots, f_k$ is a pendant edge in $H$, then $v\in D'$ (similarly as in Case 3.1) and this again contradicts the choice of $v$. Thus assume that none of the edges $f_1, \ldots, f_k$ is a pendant edge in $H$. Then, since $H$ has no good subgraph, it follows from Corollary \ref{corollary-2-small-good-subgraphs} that every vertex $v^1, \ldots, v^k$ is incident with a pendant edge in $H$. Analogously as in Case 3.1, each of the vertices  $v^1, \ldots, v^k$ belongs to $D'$ in $G$. Now, the minimality of $G$ implies in turn that the vertices $(v^i,\{v^i_{f_i}\})$ belong to $T'$ for $i \in [k]$. Consequently, the vertices $(v,\{v_{f_i}\})$ also belong to $T'$, since $T'$ is a total dominating set in $G$ and $(v,\{v_{f_i}\})$ is the only neighbor of $(v^i,\{v^i_{f_i}\})$ which is not in $D'$ (for $i \in [k]$). Finally, since all the neighbors $(v,\{v_{f_i}\})$ ($i \in [k]$) of the vertex $v$ are in $T'$, the vertex $v$ has to be in $D'$, a final contradiction proving the implication $(1) \Rightarrow (2)$.

$(2) \Rightarrow (1)$. Assume that $G=S_2(H,{\cal P},\theta)$ $($for some connected graph $H$, some family ${\cal P}=\{{\cal P}(v)\colon v\in V_H\}$ in which ${\cal P}(v)$ is a partition of the neighborhood $N_{S_2(H)}(v)$ $($for $v\in V_H$$)$, and some positive function $\theta \colon L_{S_2(H,{\cal P})}\to \mathbb{N}$$)$. If $G$ is a cycle of length 3, 6 or 9, then $G$ is a minimal DTDP-graph. Thus assume that $(V_G^o,V_G^n)$ is the only \DT-pair in $G$ and every component of $G[V_{G}^n]$ is a~star. Certainly, $G$ is a DTDP-graph (by Proposition~\ref{observ-S2(H)-jest-DPDP-grafem}), and we shall prove that $G$ is a~minimal DTDP-graph. Suppose, to the contrary, that $G$ is not a minimal DTDP-graph. Thus some proper spanning subgraph $G'$ of $G$ is a DTDP-graph. Let $e$ be any edge belonging to $G$ but not to $G'$, say $\varphi_G(e)=\{v,u\}$. Then, since $V_G^o$ is an independent set in $G$, either $|\{v,u\}\cap V_G^o|=1$ or $\{v,u\}\subseteq V_G^n$.

Let $(D',T')$ be a \DT-pair in $G'$ and, consequently, in $G-vu$ and $G$. Therefore  $(D',T')=(V_G^o,V_G^n)$ (since $(V_G^o,V_G^n)$ is the only \DT-pair in $G$) and $(V_G^o,V_G^n)$ is a \DT-pair in $G-vu$. But this is impossible, as we will see below. Assume first that $|\{v,u\}\cap V_G^o|=1$, say $v\in V_G^o$ and $u\in V_G^n$. Then $N_G(u)\cap V_G^o=\{v\}$ (by Observation \ref{pierwsze-wlasnosci-S2(H)}\,(4)) and, therefore, $N_{G-vu}(u) \cap V_G^o= \emptyset$, which contradicts the observation that  $(V_G^o,V_G^n)$ is a \DT-pair in $G-vu$. Thus assume that $\{v,u\}\subseteq V_G^n$.
Because $v$ and $u$ are adjacent in $G[V_G^n]$ and every component of  $G[V_G^n]$ is a star, at least one of the vertices $v$ and $u$ is a leaf in $G[V_G^n]$, say $v$ is a leaf in $G[V_G^n]$. Hence, $u$ is the only neighbor of $v$ in $G[V_G^n]$ and, therefore, no neighbor of $v$ belongs to $V_G^n$ in $G-vu$. Thus, $V_G^n$ is not a~total dominating set of $G-vu$, which contradicts the observation that  $(V_G^o,V_G^n)$ is a \DT-pair in $G-vu$. We conclude that $G$ is a minimal DTDP-graph.

$(1) \Rightarrow (3)$.
Assume again that $G$ is a connected minimal DTDP-graph. Then the equivalence of (1) and (2) implies that either  $G\in \{C_3, C_6, C_9\}$ or $G$ is a $2$-subdivision graph, say $G=S_2(H,{\cal P},\theta)$ (where $H$ is a connected graph of order at least 2, ${\cal P}=\{{\cal P}(v)\colon v\in V_H\}$ is a~family in which ${\cal P}(v)$ is a partition of the set  $N_{S_2(H)}(v)$ for $v\in V_H$, and $\theta \colon L_{S_2(H,{\cal P})}\to \mathbb{N}$ is a positive function). Thus, either $G=S_2(H)$ and $H\in \{C_1, C_2, C_3\}$ or $G=S_2(H,{\cal P},\theta)$, where $H$ is a connected graph of order at least 2 in which every non-pendant edge and every loop is adjacent to a pendant edge (as otherwise $H$ has a good subgraph (by Corollary \ref{corollary-small-good-subgraphs}) and then $G=S_2(H,{\cal P},\theta)$ would be a~non-minimal DTDP-graph (by Theorem \ref{zakazane-w-minimalnych-DTDP})), ${\cal P}=\{{\cal P}(v)\colon v\in V_H\}$ is a family in which every ${\cal P}(v)$ is a partition of the set  $N_{S_2(H)}(v)$ for every $v\in V_H$, and ${\cal P}$ contracts in $S_2(H,{\cal P})$ only far parts of adjacent pendant edges of $H$ (as otherwise $G=S_2(H,{\cal P},\theta)$ would be a~non-minimal DTDP-graph (by Proposition~\ref{o sklejaniu krawedzi})), and $\theta \colon L_{S_2(H,{\cal P})} \to \mathbb{N}$ is a~positive function (as we already have observed).

$(3) \Rightarrow (1)$. Assume finally that $H$, ${\cal P}$, and $\theta$ have the properties stated in (3). Then $G=S_2(H,{\cal P},\theta)$ is a DTDP-graph (by Proposition~\ref{observ-S2(H)-jest-DPDP-grafem}). We shall prove that $G$ is a minimal DTDP-graph. Since the minimality of $S_2(H,{\cal P},\theta)$ does not depend on positive values of $\theta$ (by Observation \ref{observ-S2(H)-stron-weak-supports}), we may assume that $\theta(x)=1$ for every $x\in L_{S_2(H,{\cal P})}$, and therefore we may assume that $G=S_2(H,{\cal P})$. By our assumption ${\cal P}$ contracts at most far parts of adjacent pendant edges in $H$, and since the minimality of  $S_2(H,{\cal P})$ does not depend on such contractions (by Proposition~\ref{o sklejaniu krawedzi koncowych}), we may assume that $G=S_2(H)$. We shall prove that $G=S_2(H)$ is a minimal DTDP-graph.
In order to prove this it suffices to show that $G$ has the properties stated in (2). If $H\in \{C_1, C_2, C_3\}$, then we note that $G\in \{C_3, C_6, C_9\}$. Thus, since every component of $G[V_G^n]$ is a star (in fact, a star of order 2), it remains to prove that the pair $(V_G^o,V_G^n)$ is the only \DT-pair in $G$ if every non-pendant edge and every loop is adjacent to a~pendant edge in $H$. Let $(D,T)$ be a \DT-pair in $G$. We shall prove that $(D,T)=(V_G^o,V_G^n)$. Since the pairs $(D,T)$ and $(V_G^o,V_G^n)$ form partitions of the set $V_G$, it suffices to show that $V_G^o \subseteq D$ and $V_G^n \subseteq T$. We prove these two containments showing that if $e$ is an edge (a loop, resp.) of $H$ and $\varphi_H(e)= \{v,u\}$ ($\varphi_H(e)= \{v\}$, resp.), then $\varphi_H(e) \subseteq D$ and $\{v_e,u_e\}\subseteq T$ ($\{v_e^1,v_e^2\}\subseteq T$, resp.). We distinguish three possible cases: (a) $e$ is a pendant edge in $H$; (b) $e$ joins two support vertices in $H$ (or $e$ is a~loop incident with a support vertex in $H$); (c) exactly one of the two end vertices of $e$ is a~support vertex in $H$.

{\it Case 1. The edge $e$ is a pendant edge in $H$, say $\varphi_H(e)= \{v,u\}$, where $u\in L_H$.} In this case, $(u,u_e,v_e,v)$ is a 4-path in $G=S_2(H)$ and $d_G(u_e)= d_G(v_e)=2$. Now $u\in D$ and $u_e\in T$ (by Observation \ref{observ-1}), and this implies that $v_e\in T$ (since $u_e$ belonging to a total dominating set $T$ in $G$ has a neighbor in $T$ and $v_e$ is the only neighbor of $u_e$ in $V_G\setminus D$) and $v\in D$ (since $D$ is a dominating  set in $G$ and $v$ is the only neighbor of $v_e$ which is not in $D$). Consequently, $L_H \cup S_H \subseteq D$ and, in addition, if a pendant edge $e$ joins vertices $v$ and $u$ in $H$, then $\{v_e,u_e\}\subseteq T$.

{\it Case 2. $\varphi_H(e)= \{v,u\}\subseteq S_H$.} In this case, $\varphi_H(e) \subseteq D$
(since $S_H \subseteq D$) and both $v_e$ and $u_e$ must be in $T$ (because $\{v,u\} \subseteq D$, $N_G(v_e)= \{v,u_e\}$, $N_G(u_e)=\{u,v_e\}$, and $(D,T)$ is a \DT-pair in $G$). Similarly, if $e$ is a loop incident with a~support vertex $v$ in $H$, then
$\varphi_H(e)= \{v\}\subseteq S_H \subseteq D$ and, certainly, both $v_e^1$ and $v_e^2$ are in $T$.

{\it Case 3. A non-pendant edge $e$ (which is not a loop) is incident with exactly one support vertex, say $\varphi_H(e)=\{v,u\}$ and $\varphi_H(e)\cap S_H= \{u\}$.} Let $E_H(v)$ denote the set of edges incident with $v$ in $H$. Since $v\not\in S_H$ and $e$ is a non-pendant edge, $|E_H(v)|\ge 2$ and every  element of $E_H(v)$ is a non-pendant edge (and it is not a loop). Therefore, since every non-pendant edge is adjacent to a~pendant edge in $H$, each neighbor of $v$ is a~support vertex in $H$ and, consequently, $N_H(v) \subseteq S_H \subseteq D$ in $G$. We claim that $v$ also belongs to $D$ in $G$. Suppose, to the contrary, that $v$ is in $T$. Then, since $(D,T)$ is a~\DT-pair in $G$, there is an edge $f$ in $E_H(v)$ such that $v_{f}\in D$. Suppose, without loss of generality, that $\varphi_H(f)= \{v,w\}$. Then $N_G(w_f)=\{w, v_f\}\subseteq D$, and, therefore, $N_G(w_f)\cap T=\emptyset$, which is impossible as $T$ is a total dominating set in $G$. This proves that $v\in D$ and implies that both $v$ and $u$ are in $D$. Finally, as in Case 2, we observe that both $v_e$ and $u_e$ are in $T$.
This completes the proof.  \end{proof}

As an immediate consequence of Theorem~\ref{minimal-DTDP-graphs-2}, we have the following corollaries.

\begin{cor2}
If $H$ is a graph in which every vertex is a leaf or it is adjacent to at least one leaf, then $S_2(H,{\cal P},\theta)$ is a minimal DTDP-graph if and only if ${\cal P}=\{{\cal P}(v) \colon v\in V_H\}$ is a family of partitions ${\cal P}(v)$ of sets $N_{S_2(H)}(v)$ $($$v\in V_H$$)$ which contracts only far parts of adjacent pendant edges $($if any$)$, and $\theta \colon L_{S_2(H,{\cal P})} \to \mathbb{N}$ is a~positive function.
\end{cor2}

\begin{cor2} A tree $T$ of order at least~$4$ is a minimal DTDP-graph if and only if $T$ is a~$2$-subdivision graph, say $T=S_2(R,{\cal P},\theta)$, where $R$ is a tree in which every non-pendant edge is adjacent to a~pendant edge, ${\cal P}=\{{\cal P}(v) \colon v\in V_H\}$ is a family of partitions ${\cal P}(v)$ of sets $N_{S_2(R)}(v)$ $($$v\in V_R$$)$, and ${\cal P}$ contracts only far parts of adjacent pendant edges $($if any$)$, and $\theta \colon L_{S_2(H,{\cal P})} \to \mathbb{N}$ is a positive function. \end{cor2}

\section{Open problems}

We close this paper with the following list of open problems that we have yet to settle. \\ [-24pt]
\begin{enumerate}
\item Characterize the graphs with loops which are minimal DTDP-graphs. \1
\item The {\it domatic-total domatic number} of a graph $G$, denoted ${\rm dom}_{\gamma\gamma_t}(G)$, is the maximum number of sets into which the vertex set of $G$ can be partitioned in such a way that the subgraph induced by the set is a DTDP-graph. It is clear that ${\rm dom}_{\gamma \gamma_t}(G)$ is a positive integer only for DTDP-graphs. We write  ${\rm dom}_{\gamma \gamma_t}(G)=0$ if a~graph $G$ is not a DTDP-graph. Give bounds on the domatic-total domatic number of a graph in terms of order. It is quite easy to observe that ${\rm dom}_{\gamma \gamma_t}(G)\le |V_G|/3$. For which graphs $G$ is ${\rm dom}_{\gamma \gamma_t}(G)= |V_G|/3$? If $G$ is a tree, then ${\rm dom}_{\gamma \gamma_t}(G)\le |V_G|/4$. For which trees $G$ is ${\rm dom}_{\gamma \gamma_t}(G)= |V_G|/4$. \2
\item Study relations between the set of minimal DTDP-graphs and the set of graphs $G$ for which $\gamma\gamma_t(G)=|V_G|$.
The reader interested in knowing more about the parameter $\gamma\gamma_t(G)$ is recommended to refer to the book \cite{Henning-Yeo2013}.
\end{enumerate}

\end{document}